\documentclass[11pt,a4paper]{article}
\usepackage[T1]{fontenc}
\usepackage{lmodern}
\usepackage[a4paper,margin=1in]{geometry}
\usepackage{amsmath,amssymb,amsfonts,amsthm}
\usepackage{graphicx,booktabs,array,longtable}
\usepackage{url}
\usepackage[unicode,hidelinks]{hyperref}
\usepackage[font=small,labelfont=bf,labelsep=space,hypcap=false]{caption}
\newcommand{\tablefont}{\fontsize{8.5pt}{10.5pt}\selectfont}

\newtheorem{theorem}{Theorem}
\newtheorem{lemma}[theorem]{Lemma}
\newtheorem{proposition}[theorem]{Proposition}
\theoremstyle{definition}
\newtheorem{definition}[theorem]{Definition}
\theoremstyle{plain}

\newcounter{appcertcase}

\newcommand{\PaperAuthorNames}{Yahui Yang, Erxiao Wang\textsuperscript{*}}
\newcommand{\PaperAffiliations}{Zhejiang Normal University}
\newcommand{\PaperFundingNote}{Corresponding author (\href{mailto:wang.eric@zjnu.edu.cn}{\nolinkurl{wang.eric@zjnu.edu.cn}}). This work was supported by the National Natural Science Foundation of China (NSFC--RGC grant 12361161603).}
\newcommand{\PaperFirstPageNote}{%
  \begingroup
  \renewcommand{\thefootnote}{\fnsymbol{footnote}}%
  \footnotetext[1]{\PaperFundingNote}%
  \endgroup}
\newcommand{\PaperAIStatement}{During the research and preparation of this manuscript, the authors used ChatGPT to assist with mathematical research, computational checks, and manuscript preparation, and Claude to assist with reviewing the logic of the arguments. The authors studied the relevant literature, discussed the research, and reviewed, revised, and verified the AI-assisted content. The authors take full responsibility for all content of this work, including its mathematical statements, proofs, computational results, and references.}

\begin{document}
\title{Convex monohedral tilings: a sharp neighbor bound and rigidity in the equality case}
\author{\PaperAuthorNames\\[3pt]\normalsize\PaperAffiliations}
\date{}
\maketitle
\PaperFirstPageNote

\begin{abstract}
We study both side-to-side and non-side-to-side convex monohedral tilings of the Euclidean plane in which every tile has the same number of neighbors. Reflections are allowed, point contacts count as neighbors, and no periodicity assumption is imposed. Combining boundary estimates for sums of face curvatures, angle-incidence counts and side-length matching, we prove that the common neighbor number satisfies $N\le21$ and classify all equality cases. Equality holds precisely when the prototile is a $30^\circ$-$30^\circ$-$120^\circ$ triangle and the entire tiling is similar to the central subdivision of the standard equilateral triangular tiling. This resolves the open problem on the maximum common number of neighbors recorded in Adams's \emph{The Tiling Book}.
\end{abstract}

\par

\section{Introduction}\label{sec:introduction}
Convex monohedral tilings of the Euclidean plane have been extensively studied. Their prototile classification and tile adjacency are related questions. We consider tilings in which every tile has the same number of neighbors, and determine the maximum of this common number and all tilings attaining it.

Rao's exhaustive search shows that every convex pentagon admitting a monohedral tiling of the plane belongs to one of the fifteen known families \cite{Rao2017}. Fejes T\'oth \cite[Figs. 1--9]{FejesToth1975} exhibited monohedral tilings with $6,7,8,9,10,12,14,16,21$ neighbors per tile, and asked on page 275 whether other values occur. He also described Robinson's nine-value classification under the additional assumption that all tile angles meeting at each tiling vertex are equal \cite[pp. 275--276]{FejesToth1975}. Bezdek \cite[Theorem 1]{Bezdek1980} proved that these are exactly the possible common numbers of neighbors for side-to-side convex monohedral tilings. In particular, the upper bound $21$ is known in that setting. For convex monohedral tilings, he also obtained the corresponding classification when each tile boundary contains the same number of tiling vertices \cite[Theorem 2]{Bezdek1980}. Congruent tiles may have different subdivisions of their sides, so neither the side-to-side condition nor the condition on boundary vertex counts follows from congruence. Allowing two prototiles, Friedman and van der Zwaag \cite{FriedmanZwaag2002} constructed tilings with $15$, $32$ and $43$ neighbors per tile; the last two use nonconvex tiles.

Neighbor numbers have also been studied for packings, which allow gaps between tiles. Early work includes Fejes T\'oth \cite{FejesToth1969}; see also the surveys \cite[Section 2.5]{FejesToth2018} and \cite[Section 13.6]{Neighbors2022}. Linhart obtained bounds depending on the maximum number of disks meeting at a point, for finite packings and for locally finite infinite packings in which the ratio of boundary disks to interior disks tends to zero \cite[S\"atze 1--2]{Linhart1977}. The latter conditions hold for the tilings considered here, but these estimates alone do not yield the uniform bound $21$. The Newton number $\nu(P)$ of a convex disk $P$ is the largest number of congruent copies that can touch one copy in a packing by congruent copies of $P$, with pairwise disjoint interiors; point contacts are included. Such a packing is maximal in the Newton-number sense if every member has $\nu(P)$ neighbors. G\'acs \cite{Gacs1972} proved an absolute upper bound for the number of neighbors in such packings. For Newton-maximal packings, Linhart \cite[Satz 3]{Linhart1977} proved a uniform upper bound of $21$ neighbors per disk. Wegner \cite{Wegner1992} proved that the Newton number of the $30^\circ$-$30^\circ$-$120^\circ$ triangle is exactly $21$. Thus the known $21$-neighbor tiling is also maximal in the Newton-number sense.

In Section 2.2, page 97 of \emph{The Tiling Book}, Adams lists the maximum common number of neighbors in a convex monohedral tiling, without a side-to-side assumption, as an open question \cite{Adams2022}. The problem does not require the common neighbor number to attain the Newton number of the prototile.

We consider tilings of the plane by congruent copies of one convex prototile, allowing reflections. The copies cover the plane and have pairwise disjoint interiors; each copy is called a tile. Two distinct tiles are neighbors if their intersection is nonempty, including a single point, as in \cite[Definition 1.7]{Adams2022}. If each tile has exactly $N$ neighbors, we call $N$ the common number of neighbors. Our setting includes non-side-to-side tilings, in which a corner of one tile may lie in the relative interior of a side of another. These subdivision vertices do not change the number of sides of the prototile. A tiling vertex is a point incident to at least three tiles, and its degree is the number of incident tiles. It is well known that a convex polygon admitting a monohedral tiling of the plane must be a triangle, quadrilateral, pentagon or hexagon \cite[Theorem 2.18, p.~151]{Adams2022}.

\par\addvspace{6.5bp}\noindent\begin{minipage}{\linewidth}
Our main result is the following.\par\smallskip
\begin{theorem}[Neighbor bound and equality case]\label{thm:main}
Let $P\subset\mathbb R^2$ be a compact convex set with nonempty interior whose congruent copies tile the plane with pairwise disjoint interiors. Reflections are allowed. If every tile has exactly $N$ neighbors, then $N\le21$. Equality $N=21$ holds if and only if $P$ is a $30^\circ$-$30^\circ$-$120^\circ$ triangle and the entire tiling is similar to the central subdivision of the standard equilateral triangular tiling.
\end{theorem}
\end{minipage}\par\addvspace{6.5bp}

The central subdivision in the theorem is obtained by joining the center of each equilateral triangle to its three vertices, producing three congruent $30^\circ$-$30^\circ$-$120^\circ$ triangles. This construction appears in Fejes T\'oth \cite[Fig. 9]{FejesToth1975}. As shown in Fig.~\ref{fig:center-trisection}, each tile has $3$ neighbors sharing a segment and $18$ neighbors meeting it only at a vertex, so its common number of neighbors is $21$.

\par\addvspace{\intextsep}
\noindent\begin{minipage}{\linewidth}
\makeatletter\def\@captype{figure}\makeatother
\centering
\includegraphics[width=95mm]{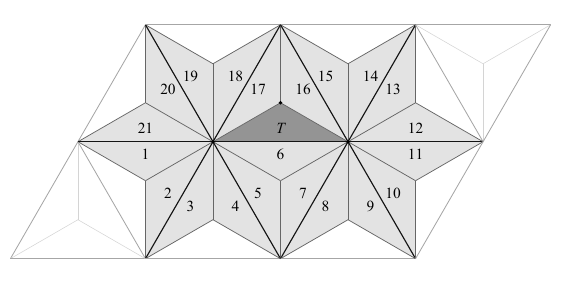}
\caption{Central subdivision of the standard equilateral triangular tiling}\label{fig:center-trisection}\label{fig:construction}
\end{minipage}\par\addvspace{\intextsep}

Theorem~\ref{thm:main} resolves the open problem recorded in Adams's \emph{The Tiling Book} and determines both the prototile and the global arrangement in the equality case. The complete set of realizable common neighbor numbers in the interval $6\le N\le21$, together with corresponding constructions, will be studied in future work.

Face curvature, estimates for reciprocal sums under degree constraints, and boundary estimates for finite regions have been used in Linhart's study of packings \cite[pp.~287--290]{Linhart1977}. We combine these tools with discharging, angle-incidence counts \cite[proof of Proposition 3]{Rao2017}, counting constraints induced by unmatched side lengths, and local configuration analysis to prove the upper bound on the common number of neighbors. The local configurations in the equality case and their propagation determine the entire tiling.

Section~\ref{sec:basic-facts} establishes the neighbor-count identity, the curvature boundary estimate, limiting weights and the injective count. Section~\ref{sec:main-proof} proves the main theorem. Section~\ref{sec:polygons} treats convex quadrilaterals, pentagons and hexagons. Section~\ref{sec:triangles} treats triangles, first according to whether the smallest angle is below $30^\circ$, and then by considering the nonexceptional parameters in the four families, the $22$ ordinary exceptions and the two special exceptions. Section~\ref{sec:rigidity} classifies the equality case $N=21$. The appendix following the references gives the exact slacks in the finite inequalities and the successive vertex-type elimination data.

\section{Basic facts}\label{sec:basic-facts}

Throughout the paper, we work under the assumptions of Theorem~\ref{thm:main}. Let $A>0$ and $D$ be the area and diameter of the prototile. Every tile meeting a bounded set $K$ lies in its $D$-neighborhood. Area comparison therefore shows that only finitely many tiles meet $K$, so the tiling is locally finite. The boundary of each tile is covered by its intersections with finitely many neighbors. By convexity, each such intersection is a segment or a point, and hence the prototile is a convex polygon. By the side-number restriction recalled in the introduction, the prototile is a convex $n$-gon with $3\le n\le6$. Subdivision vertices in the relative interiors of sides are counted separately and do not change the number of sides of the prototile.

We use \emph{corners} and \emph{sides} for a tile, and \emph{vertices} and \emph{edges} for the tiling, following \cite[Section 2]{LiangChenWang2026}. A point belonging to at least three tiles is called a \emph{tiling vertex}. Its \emph{degree} is the number of incident tiles, denoted by $q(v)\ge3$. The angles of these tiles at $v$ sum to $2\pi$, and at most one tile contributes a straight angle $\pi$. A tiling vertex with no straight angle is called a \emph{full vertex}; one with a straight angle is called a \emph{half-vertex}, and the tile contributing this straight angle is called its \emph{straight-angle tile}. A half-vertex lies in the relative interior of a side of its straight-angle tile and is a corner of every other incident tile. Every corner of a tile is a tiling vertex. A tiling vertex in the relative interior of a side of a tile is called a \emph{subdivision vertex} of that tile. If the intersection of two tiles is a nondegenerate segment, its relative interior contains no tiling vertex. A segment between consecutive tiling vertices on the boundary of a tile is called an \emph{edge} of the tiling. Every edge belongs to exactly two tiles, and the boundary of each tile is a finite simple closed chain of these edges. A tiling is \emph{side-to-side}, also called edge-to-edge, if each tile side is one edge of the tiling; otherwise it is \emph{non-side-to-side}.

\par\addvspace{6.5bp}\noindent\begin{minipage}{\linewidth}
\begin{lemma}[Neighbor-count identity]\label{lem:neighbors}
For a tile $T$, let $\mathcal V(T)$ be the set of tiling vertices on its boundary, and put $m(T)=|\mathcal V(T)|$. Then
\[
N=\sum_{v\in\mathcal V(T)}(q(v)-2). \tag{1}
\]
\end{lemma}
\end{minipage}\par

\begin{proof}
If a contact point belonged to only two tiles, local finiteness and covering would force both to have a straight angle there and to share a segment extending on both sides. Thus every isolated contact and every endpoint of a shared segment is a tiling vertex. At each $v\in\mathcal V(T)$, there are $q(v)-1$ tiles other than $T$. Summing over these vertices gives $\sum_{v\in\mathcal V(T)}(q(v)-1)$ incidences with neighbors: a point-contact neighbor is counted once, and a segment-sharing neighbor is counted at both endpoints. By convexity, the $m(T)$ tiling edges of $T$ correspond to $m(T)$ distinct segment-sharing neighbors. Subtracting these $m(T)$ repeated incidences gives equation~(1). Since $q(v)\ge3$, equation~(1) also gives $m(T)\le N$ and $q(v)\le N$.
\end{proof}

Let $B_R$ be the closed disk of radius $R$ centered at a fixed origin, let $\mathcal F_R$ be the set of tiles contained in $B_R$, and put $f_R=|\mathcal F_R|$. These tiles cover $B_{R-D}$. Area comparison gives
\[
\pi(R-D)^2\le Af_R\le\pi R^2\qquad(R\ge D),
\qquad f_R=\pi R^2/A+O(R).
\]
Tiles meeting an annulus of fixed width lie in another annulus of fixed width, so area comparison gives $O(R)$ such tiles. Each tile has at most $N$ boundary vertices, so these tiles are incident to only $O(R)$ vertices. A counting correspondence or curvature transfer therefore produces an error of $O(R)$ across the disk boundary whenever its distances, multiplicities and weights are uniformly bounded. The constants may depend on the prototile, $N$, and these uniform bounds.

For Euler counting in non-side-to-side spherical triangle tilings, see \cite[Section 2, equations (2.1)--(2.2)]{LiangChenWang2026}; below we carry out analogous counting in finite planar regions and estimate the boundary terms.

\par\addvspace{6.5bp}\noindent\begin{minipage}{\linewidth}
\begin{lemma}[Estimate for the sum of face curvatures]\label{lem:curvature}
For a tile $T$, define its face curvature by
\[
\kappa(T)=1-\frac{m(T)}2+\sum_{v\in\mathcal V(T)}\frac1{q(v)}.
\]
There are constants $C>0$ and $R_0>0$, depending only on the fixed prototile and $N$, such that
\[
\left|\sum_{T\in\mathcal F_R}\kappa(T)\right|\le CR
\qquad(R\ge R_0). \tag{2}
\]
The same estimate holds with $\kappa$ replaced by the redistributed curvature $\widetilde\kappa$, provided that each transfer subtracts an amount from one tile and adds it to another, and the Euclidean distances between the two tiles, the absolute amounts transferred, and the total number of transfers involving any one tile are uniformly bounded. The constants may then also depend on these bounds.
\end{lemma}
\end{minipage}\par

\begin{proof}
Let $\Gamma_R$ be the finite plane graph formed by the boundaries of the tiles in $\mathcal F_R$, and let $v_R,e_R$ be its numbers of vertices and edges. Every unselected tile has an interior point outside $B_R$, and its connected interior does not cross $\Gamma_R$. Thus $\Gamma_R$ has exactly $f_R$ bounded faces, one for each selected tile. Edges and vertices incident to both selected and unselected tiles lie in an annulus of fixed width and number $O(R)$. Each connected component contains such a boundary vertex, so the number of components is also $O(R)$. Euler's formula gives
\[
v_R-e_R+f_R=O(R).
\]
Away from these boundary terms, every edge is incident to two selected tiles and every vertex has all its incident tiles selected. Hence
\[
\sum_{T\in\mathcal F_R}m(T)=2e_R+O(R),\qquad
\sum_{T\in\mathcal F_R}\sum_{v\in\mathcal V(T)}\frac1{q(v)}
=v_R+O(R).
\]
Substitution into the curvature definition yields
\[
\sum_{T\in\mathcal F_R}\kappa(T)
=f_R-e_R+v_R+O(R)=O(R),
\]
which proves inequality~(2). Since $f_R=\pi R^2/A+O(R)$, we also obtain
\[
\lim_{R\to\infty}\frac1{f_R}\sum_{T\in\mathcal F_R}\kappa(T)=0.
\]
Transfers between selected tiles cancel. Transfers across the disk boundary involve only $O(R)$ selected tiles in an annulus of fixed width. The number of transfers per tile and their amounts are uniformly bounded, so the total boundary contribution is $O(R)$. Thus $\widetilde\kappa$ also satisfies inequality~(2), and its average tends to zero.
\end{proof}

\par\addvspace{6.5bp}\noindent\begin{minipage}{\linewidth}
\begin{lemma}[Local curvature bound]\label{lem:local-curvature}
For a convex polygonal tile $T$, put $m=m(T)\ge n$. Then
\[
\kappa(T)\le\frac{4-m}{6}+\frac1{N-m+3}\qquad(3\le m\le N). \tag{3}
\]
For fixed $N$, the right-hand side of inequality~(3) is strictly decreasing in $m$.
\end{lemma}
\end{minipage}\par

\begin{proof}
The neighbor-count identity~(1) gives
\[
\sum_{v\in\mathcal V(T)}q(v)=N+2m\qquad(q(v)\ge3).
\]
For any $x,y\ge3$, the inequality $(x-3)(y-3)\ge0$ gives
\[
\frac1x+\frac1y=\frac{x+y}{xy}
\le\frac{x+y}{3(x+y-3)}
=\frac13+\frac1{x+y-3}.
\]
Thus replacing two entries by $3$ and their sum minus $3$ preserves the total and does not decrease the reciprocal sum. Repeating this replacement leaves $m-1$ entries equal to $3$. Hence
\[
\sum_{v\in\mathcal V(T)}\frac1{q(v)}
\le\frac{m-1}{3}+\frac1{N+2m-3(m-1)}
=\frac{m-1}{3}+\frac1{N-m+3}.
\]
Substitution into the curvature definition in Lemma~\ref{lem:curvature} gives inequality~(3).
For $3\le m\le N$, the monotonicity follows from
\[
\frac{\mathrm d}{\mathrm dm}
\left(\frac{4-m}{6}+\frac1{N-m+3}\right)
=-\frac16+\frac1{(N-m+3)^2}\le-\frac1{18}<0.
\]
\end{proof}

The following notation and counting tools will be used for triangular prototiles in Sections~\ref{sec:triangles} and~\ref{sec:rigidity}.

Let $\alpha\le\beta\le\gamma$ be the angles of the triangular prototile, measured in degrees. Fix labels for its three corners and transfer them to each tile through a chosen congruence with the prototile. We retain these labels even when some angle values are equal. Let $r,s,t$ be the multiplicities of the three labeled angles at a vertex, and put
\[
\varepsilon=\begin{cases}
0,&\text{at a full vertex, with no straight angle},\\
1,&\text{at a half-vertex, with one straight angle}.
\end{cases}
\]
The non-straight angles at the vertex sum to $180^\circ(2-\varepsilon)$. A quadruple $\tau=(r,s,t;\varepsilon)$ satisfying
\[
r,s,t\in\mathbb Z_{\ge0},\quad \varepsilon\in\{0,1\},\quad
r\alpha+s\beta+t\gamma=180^\circ(2-\varepsilon),\quad q=r+s+t+\varepsilon\ge3
\]
is called a \emph{candidate vertex type}. Since $\alpha>0$ and $r+s+t\le360^\circ/\alpha$, there are finitely many candidate vertex types, some of which may be absent from the tiling.

Let $v_\tau(B_R)$ be the number of tiling vertices of type $\tau$ in $B_R$. Every vertex is a corner of some tile, and every tile meeting $B_R$ lies in $B_{R+D}$. Hence
\[
\sum_\tau v_\tau(B_R)\le3f_{R+D}.
\]
The tile-count estimate bounds $x_\tau(R)=v_\tau(B_R)/f_R$ for sufficiently large $R$. Since there are finitely many candidate vertex types, we can take simultaneous limits along a common subsequence of radii $R_j\to\infty$. Write
\[
x_\tau=\lim_{j\to\infty}x_\tau(R_j)\ge0,\qquad
h=\sum_\tau \varepsilon x_\tau.
\]
Fix this subsequence from now on. We call $x_\tau$ the \emph{limiting weight} of the vertex type $\tau$, and $h$ is the limiting number of half-vertices normalized by the number of tiles. If $x_\tau>0$, this vertex type must actually occur, but an occurring vertex type may have zero limiting weight. All subsequent sums over vertex types run over all candidate vertex types, with all quantities evaluated at the corresponding vertex type.

\par\addvspace{6.5bp}\noindent\begin{minipage}{\linewidth}
\begin{lemma}[Side-length mismatch and injective counting]\label{lem:unmatched}
Suppose the prototile is a scalene triangle with side lengths $a<b<c$ opposite $\alpha,\beta,\gamma$, respectively. If the angle multiplicities $r,s,t$ at a full vertex $v$ do not all have the same parity, then some ray from $v$ separates two incident tiles whose sides along the ray have different lengths. For each such $v$, choose one such ray and let $v'=\varphi(v)$ be the other endpoint of the shorter side. Then $v'$ is a half-vertex, and every half-vertex has at most two preimages; if no side length is twice another, then $\varphi$ is injective.
\end{lemma}
\end{minipage}\par

\begin{proof}
At $v$, an $\alpha$-angle is bounded by sides of lengths $b,c$, a $\beta$-angle by sides of lengths $a,c$, and a $\gamma$-angle by sides of lengths $a,b$. Thus the numbers of incident side ends of the three lengths are
\[
s+t,\qquad r+t,\qquad r+s.
\]
If $r,s,t$ do not all have the same parity, the number of incident side ends is odd for at least one side length. Equal-length sides therefore cannot all be paired, so some ray has different side lengths on its two sides. The other endpoint $v'$ of the shorter side lies strictly inside the longer side, whose tile provides a straight angle there. Thus $v'$ is a half-vertex.

Every half-vertex $v'$ lies in the relative interior of a side of the unique tile contributing the straight angle at $v'$. Since the preimages are full vertices, they can only be the two endpoints of that side. Thus there are at most two preimages. If both endpoints $v_1,v_2$ map to $v'$, let $L$ be the length of the side, and let $l_1,l_2$ be the distances from its endpoints to $v'$. Then
\[
L=l_1+l_2,\qquad l_1,l_2<L,\qquad L,l_1,l_2\in\{a,b,c\}.
\]
If $l_1\ne l_2$, then $L,l_1,l_2$ are the three distinct side lengths of the prototile, contradicting the triangle inequality $L<l_1+l_2$. Hence $l_1=l_2$, so $L=2l_1$. Thus $\varphi$ is injective when no side length is twice another.

The construction gives $|v-v'|\le D$ for a vertex $v$ and its image $v'=\varphi(v)$. Thus counting these full vertices and their images in a disk incurs only an $O(R)$ boundary error. Since every half-vertex has at most two preimages, normalization and passage to the common subsequence show that the total weight of these full vertices is at most $2h$.
\end{proof}

\noindent\textbf{Remark.} Suppose that the angles of a triangle are rational multiples of $\pi$ and that one side is twice another. Write $\theta$ for the angle opposite the longer side and $\psi$ for the angle opposite the shorter side, measured in radians. The sine rule gives $\sin\theta=2\sin\psi$. Set
\[
\xi=\sin\theta=2\sin\psi=-i(e^{i\psi}-e^{-i\psi}),
\]
and take a cyclotomic field $K$ containing $e^{i\psi},e^{i\theta},i$. Since $e^{i\psi},e^{i\theta},i$ are roots of unity, $\xi$ is an algebraic integer. Every embedding of $K$ maps roots of unity to roots of unity, so every image of $\xi$ has the form $\pm\sin(k\theta)$ and hence has absolute value at most $1$. Since $\theta$ is an interior angle of a triangle, $0<\xi\le1$.

If $0<\xi<1$, the factor $\xi$ corresponding to the identity embedding has absolute value less than $1$, while all other factors have absolute value at most $1$. The field norm of $\xi$ therefore has absolute value strictly less than $1$. This contradicts the nonzero integer norm of the nonzero algebraic integer $\xi$ \cite[Theorem 2.1 and Corollaries 2.20--2.21]{Milne2020}. Hence $\xi=1$, so $\theta=\pi/2$. The conditions $\sin\psi=1/2$ and $\psi+\theta<\pi$ then give $\psi=\pi/6$. Therefore the only rational-angle triangle with one side twice another is the $30^\circ$-$60^\circ$-$90^\circ$ triangle. Every scalene rational-angle triangle other than this right triangle satisfies the injectivity condition in Lemma~\ref{lem:unmatched}. The following parity correction records this side-length restriction.

\par\addvspace{6.5bp}\noindent\begin{minipage}{\linewidth}
\begin{definition}[Parity correction]\label{def:parity}
For a vertex type $\tau=(r,s,t;\varepsilon)$, define
\[
P(\tau)=
\begin{cases}
1,&\varepsilon=0,\quad r,s,t\text{ do not all have the same parity},\\
0,&\varepsilon=0,\quad r,s,t\text{ are all odd or all even},\\
-1,&\varepsilon=1.
\end{cases}
\]
\end{definition}
\end{minipage}\par

When the prototile is scalene and no side is twice another, Lemma~\ref{lem:unmatched} gives an injection from the full vertices with $P=1$ into the half-vertices with $P=-1$, with correspondence distance at most $D$. In comparing the two counts in $B_R$, only correspondences within distance $D$ of the boundary can cross it, and there are $O(R)$ of these. Normalizing and passing to the chosen subsequence gives
\[
\sum_\tau P(\tau)x_\tau\le0. \tag{4}\label{eq:parity-bound}
\]

\section{Proof of the main theorem}\label{sec:main-proof}

\subsection{Convex quadrilaterals, pentagons and hexagons}\label{sec:polygons}

\par\addvspace{6.5bp}\noindent\begin{minipage}{\linewidth}
\begin{proposition}\label{prop:polygons}
If the prototile is a convex pentagon or hexagon, then $N\le8$. If it is a convex quadrilateral, then $N\le14$.
\end{proposition}
\end{minipage}\par

\begin{proof}
Fix a tile $T$ and write $m=m(T)$. If $n\in\{5,6\}$ and $N\ge9$, then $m\ge n\ge5$, so Lemma~\ref{lem:local-curvature} gives $\kappa(T)\le-1/6+1/(N-2)\le-1/42$ for every tile, contrary to inequality~(2) in Lemma~\ref{lem:curvature}. Thus convex pentagons and hexagons both satisfy $N\le8$.

For convex quadrilaterals, suppose that $N\ge15$ and put $H=N-1$. When $m=4$, a tile of type $(H,3,3,3)$ has curvature $1/H$. Every other degree list with $m=4$ has at least two entries at least $4$, so its three smallest entries are bounded below by $3,3,4$, respectively. Equation~(1) gives a sum of $N+8$ for all four degrees. Applying the reciprocal-sum estimate in the proof of Lemma~\ref{lem:local-curvature} (cf.\ the related counting in \cite[Lemma 2.1]{LiaoQianWangXu2024}) gives
\[
\sum_{v\in\mathcal V(T)}\frac1{q(v)}
\le\frac23+\frac14+\frac1{N-2},\qquad
\kappa(T)\le-\frac1{12}+\frac1{N-2}<0.
\]
For $m\ge5$, the curvature bound in Lemma~\ref{lem:local-curvature} decreases with $m$ and is negative at $m=5$, so the curvature is negative. Thus only tiles of type $(H,3,3,3)$ have positive curvature.

Let $v$ be a vertex of degree $H$. For any tile $T$ incident to $v$, each of its other boundary vertices contributes at least $1$ to the sum in Lemma~\ref{lem:neighbors}. Hence
\[
N\ge(H-2)+(m(T)-1)=N+m(T)-4.
\]
Together with $m(T)\ge4$, this gives $m(T)=4$. The other three vertices all have degree $3$, and no side has a subdivision vertex. Thus all $H$ tiles around $v$ have positive curvature. In cyclic order around $v$, write them as
\[
Q_i=(v,a_i,b_i,a_{i+1}),\qquad i\pmod H.
\]
Fig.~\ref{fig:sphere-gluing}(a) shows the local configuration of the inner and outer tiles. Convexity and the cyclic order of the inner tiles imply that all $a_i,b_i$ are distinct. Let $G_i$ be the third tile at $a_i$; its boundary contains the path $(b_{i-1},a_i,b_i)$. If $G_i=G_{i+1}$, this tile shares two noncollinear sides with $Q_i$, forcing the interior overlap shown in the figure. Thus $G_i\ne G_{i+1}$. Only two inner tiles meet at $a_i$, so no $G_i$ is an inner tile.

If $m(G_i)=4$, Lemma~\ref{lem:neighbors} gives $q(v')=H$ at its fourth vertex $v'$. The inner tiles fill a neighborhood of $v$, and their other vertices have degree $3$, so $v'$ is distinct from these vertices. Since $q(b_i)=3$ and $G_i$ has no subdivision vertices, $G_i$ and $G_{i+1}$ share the entire side $b_iv'$. As $q(v')=H$, $G_{i+1}$ also has positive curvature and no subdivision vertices. Continuing around the boundary gives $G_j=(v',b_{j-1},a_j,b_j)$. These outer tiles are distinct. The two layers are disks and glue along their entire common boundary to form a sphere (Fig.~\ref{fig:sphere-gluing}(b)), contradicting their planar embedding. Hence every outer tile satisfies $m(G_i)\ge5$.

\begin{figure}[htbp]
\centering
\includegraphics[width=120mm]{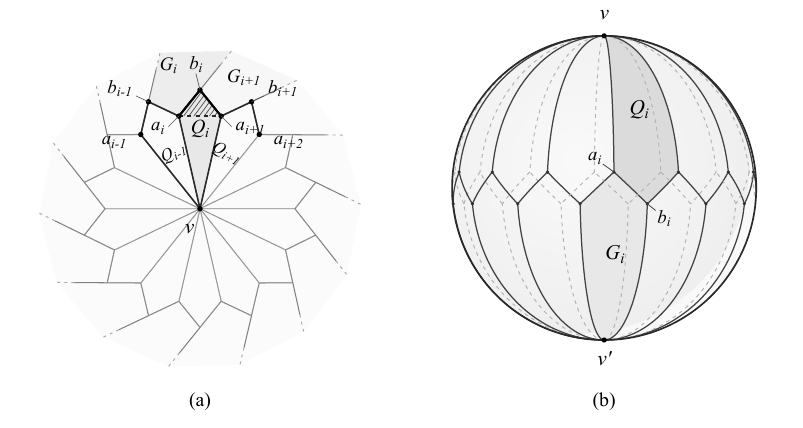}
\caption{Local configuration and gluing to a sphere}\label{fig:sphere-gluing}
\end{figure}

At every vertex of degree $H$, transfer the curvature $1/H$ of each $Q_i$ to $G_i$ as in Fig.~\ref{fig:sphere-gluing}(a), with each transfer corresponding to the path $(b_{i-1},a_i,b_i)$ on the receiving tile's boundary. Within one star, since $G_i\ne G_{i+1}$, the indices corresponding to the same receiving tile are pairwise nonadjacent, so their paths have no common vertex. Since $q(a_i)=q(b_i)=3$ and no outer tile has positive curvature, every positive-curvature tile through these vertices belongs to the current star. Paths from distinct stars therefore also have no common vertex. Thus a tile with $m$ boundary vertices receives at most $\lfloor m/3\rfloor$ transfers.

After discharging, every originally positive-curvature tile has curvature zero. Each receiving tile satisfies $m\ge5$. Lemma~\ref{lem:local-curvature} and the bound on incoming transfers give the following upper bound for its redistributed curvature:
\[
\frac{4-m}{6}
+\frac1{N-m+3}
+\frac{\lfloor m/3\rfloor}{N-1}.
\]
The right-hand side is negative at $m=5$. For $m\ge6$, replace $\lfloor m/3\rfloor$ by $m/3$. The derivative estimate in the proof of Lemma~\ref{lem:local-curvature}, together with $N\ge15$, shows that this upper bound decreases with $m\le N$; it is negative at $m=6$. Thus the curvature of each receiving tile remains negative. The curvatures of the other negative-curvature tiles, which receive no transfers, are unchanged. Thus tiles outside the positive-curvature stars all retain negative curvature, with a uniform negative upper bound.

Let $f_R^{+}$ count the tiles in $\mathcal F_R$ with positive curvature before discharging, and let $f_R^{-}$ count the remaining tiles. Each positive-curvature tile sends one transfer, and each receiving tile accepts at most $\lfloor N/3\rfloor$ transfers. Transfers across the disk boundary contribute only $O(R)$, so
\[
f_R=f_R^{+}+f_R^{-}\le(1+\lfloor N/3\rfloor)f_R^{-}+O(R).
\]
The preceding inequality and the tile-count estimate in Section~\ref{sec:basic-facts} give a uniform positive lower bound for $f_R^{-}/f_R$ for sufficiently large $R$. Consequently, the average curvature after discharging has a uniform negative upper bound, contrary to inequality~(2) in Lemma~\ref{lem:curvature}. Hence $N\le14$.
\end{proof}

\subsection{The neighbor bound for triangles}\label{sec:triangles}\label{sec:finite-classification}\label{sec:finite-inequalities}

For triangular prototiles, we retain the angle labels, candidate vertex types and common subsequence fixed in Section~\ref{sec:basic-facts}.

\begin{lemma}[Angle-incidence counts for triangles]\label{lem:angle-incidence}
The limiting weights satisfy the following two sets of identities.
\par\noindent\textup{(i)}~The angle-incidence counts satisfy
\[
\sum_\tau r x_\tau=\sum_\tau s x_\tau=\sum_\tau t x_\tau=\sum_\tau (2-\varepsilon) x_\tau=1.
\]
\par\noindent\textup{(ii)}~The vertex degrees and the neighbor number satisfy
\[
\sum_\tau(q-2)x_\tau=2,\qquad
\sum_\tau(q-3)x_\tau=\frac{3-h}{2},\qquad
N=\sum_\tau q(q-2)x_\tau. \tag{5}\label{eq:angle-incidence}
\]
\end{lemma}

\begin{proof}
Each tile contributes once at each of its three labeled corners, and its non-straight angles sum to $180^\circ$. Counting incidences in $B_R$ by tiles and by vertices gives a difference supported in an annulus of fixed width, and hence of order $O(R)$. Dividing by $f_R$ and passing to the fixed subsequence gives the first three angle-incidence identities. Counting the total non-straight angle in the same way gives $\sum_\tau (2-\varepsilon) x_\tau=1$. Since $q=r+s+t+\varepsilon$, we have
\[
\sum_\tau qx_\tau=3+h,\qquad 2\sum_\tau x_\tau=1+h,
\]
which together give the first two identities in equation~(\ref{eq:angle-incidence}). Next, sum equation~(1) in Lemma~\ref{lem:neighbors} over the tiles and interchange the order of summation. A vertex of degree $q$ is counted once by each of its $q$ incident tiles, contributing $q-2$ each time. Hence
\[
N f_R=\sum_{v\in B_R}q(v)(q(v)-2)+O(R).
\]
Dividing by $f_R$ and passing to the fixed subsequence gives
\[
N=\sum_\tau q_\tau(q_\tau-2)x_\tau.
\]
\end{proof}

Suppose that every vertex type of positive weight has degree at most $d$. Summing $q(q-2)=3(q-2)+(q-3)(q-2)$ with the limiting weights and using~(\ref{eq:angle-incidence}) and $3\le q\le d$ gives
\[
N\le6+\frac{(d-2)(3-h)}2\le\frac{3d}2+3. \tag{6}\label{eq:degree-bound}
\]
If $\alpha\ge30^\circ$, then every vertex has $q\le12$, so $N\le21$. We therefore need only consider $\alpha<30^\circ$.

\par\addvspace{6.5bp}\noindent
\begin{minipage}{\linewidth}
\begin{lemma}[Discharging at half-vertices of degree four]
\label{lem:low-degree}
If every vertex satisfies $q\ge4$ and every vertex of degree four is a half-vertex, then $N\le14$.
\end{lemma}
\end{minipage}\par

\begin{proof}
Suppose, to the contrary, that $N\ge15$. At each half-vertex of degree four, transfer curvature $1/24$ from each of the three tiles having a corner there to the tile providing the straight angle, which thus receives a total of $1/8$ at that vertex. If a tile $T$ has no subdivision vertices, then $m(T)=3$. Let $j$ be the number of its corners of degree four. By equation~(1), the three corner degrees sum to $N+6$, so $j\le2$. For $j=2,1,0$, maximizing the reciprocal sum as in the proof of Lemma~\ref{lem:local-curvature} gives the following respective upper bounds on the curvature after discharging:
\[
\frac1{13}-\frac2{24}=-\frac1{156},\quad
-\frac12+\frac14+\frac15+\frac1{12}-\frac1{24}
=-\frac1{120},\quad
-\frac12+\frac25+\frac1{11}=-\frac1{110}.
\]

If $T$ has subdivision vertices, put $k=m(T)-3\ge1$. The same reciprocal-sum argument, now with $q\ge4$, bounds its original curvature by $-k/4+1/(N-2k-2)$. Each subdivision vertex allows $T$ to receive at most $1/8$, so its curvature after discharging is at most $-k/8+1/(N-2k-2)$. Using $N\ge15$ and $N\ge2(3+k)$, for $k=1,2$ this bound is at most $-3/88,-5/36$, respectively, and for $k\ge3$ it is at most $-1/8$. Thus all tile curvatures after discharging have a common negative upper bound, contrary to inequality~(2).
\end{proof}

By Lemma~\ref{lem:low-degree}, $N>14$ forces a vertex of degree three or a full vertex of degree four. At a full vertex of degree three, if the three angle values are distinct, their sum is $180^\circ$; if one value $x$ occurs twice and the remaining value is $y\ne x$, then $2x+y<360^\circ$. Thus all three angles must be equal to $120^\circ$. At a half-vertex of degree three, the other two angles sum to $180^\circ$. Since two distinct angle values of the prototile sum to less than $180^\circ$, both angles must be $90^\circ$.

At a full vertex of degree four, if two distinct values $x,y$ each occur twice, then $2(x+y)<360^\circ$. If all three distinct values occur, the angle sum is $180^\circ$ plus one interior angle, again less than $360^\circ$. If the four angles are equal, each is $90^\circ$. The remaining case is $3x+y=360^\circ$, where $x\ne y$. If $x$ is not the largest angle, then $x\le\beta<90^\circ$ and $x+y<180^\circ$, giving $3x+y<360^\circ$, a contradiction. Hence $x=\gamma$, leaving only $3\gamma+\alpha=360^\circ$ or $3\gamma+\beta=360^\circ$. This classification by angle value also includes isosceles triangles; the incidence counting still distinguishes the three corner positions. Combining the restrictions from these two kinds of low-degree vertices, when $N>14$ and $\alpha<30^\circ$ we must have $\gamma=90^\circ$, $3\gamma+\alpha=360^\circ$, $3\gamma+\beta=360^\circ$ or $\gamma=120^\circ$. These give families R, A, B and C, respectively, in Table~\ref{tab:1}.

\begingroup
\tablefont
\setlength{\tabcolsep}{3pt}
\renewcommand{\arraystretch}{1}
\setlength{\LTcapwidth}{\textwidth}

\begin{longtable}{
>{\centering\arraybackslash}p{\dimexpr 0.3000\textwidth-2\tabcolsep\relax}
>{\centering\arraybackslash}p{\dimexpr 0.7000\textwidth-2\tabcolsep\relax}}
\caption{Four families obtained from vertices of degrees three and four}
\label{tab:1}\\
\toprule
Family & Angles, where $0<\alpha<30^\circ$ \\
\midrule
\endfirsthead
\caption[]{\textit{Continued}}\\
\toprule
Family & Angles, where $0<\alpha<30^\circ$ \\
\midrule
\endhead
\midrule
\multicolumn{2}{r}{\textit{Continued on next page}}\\
\endfoot
\bottomrule
\endlastfoot
R (right triangles) & $(\alpha,90^\circ-\alpha,90^\circ)$ \\
A & $(\alpha,60^\circ-2\alpha/3,120^\circ-\alpha/3)$ \\
B & $(\alpha,90^\circ-3\alpha/2,90^\circ+\alpha/2)$ \\
C & $(\alpha,60^\circ-\alpha,120^\circ)$ \\
\end{longtable}
\endgroup

We use the balance $r-s$ for family R, $3(2-\varepsilon)-s-2t$ for families A and C, and $2(2-\varepsilon)-s-t$ for family B. By Lemma~\ref{lem:angle-incidence}\textup{(i)}, each balance has weighted sum zero. As in \cite[proof of Proposition 3]{Rao2017}, a balance that is nonnegative on every candidate vertex type must therefore vanish on every vertex type of positive limiting weight.

A parameter is \emph{exceptional} if it admits a candidate vertex type of negative balance. Combining the negative-balance condition with the vertex angle equation for each family and using $0<\alpha<30^\circ$ gives all integer solutions listed in Tables~\ref{tab:2}--\ref{tab:5}.

For family R, the vertex angle equation is
\[
(r-s)\alpha=90^\circ(2(2-\varepsilon)-s-t).
\]
If $r<s$, the preceding equation gives $s+t-2(2-\varepsilon)>0$, and
\[
\alpha=\frac{90^\circ(s+t-2(2-\varepsilon))}{s-r}<30^\circ
\quad\Longleftrightarrow\quad r+2s+3t<6(2-\varepsilon).
\]
There are no solutions when $\varepsilon=1$. All integer solutions for $\varepsilon=0$ are listed in Table~\ref{tab:2}.

\begingroup
\tablefont
\setlength{\tabcolsep}{2pt}
\renewcommand{\arraystretch}{1}
\setlength{\LTpre}{2pt}
\setlength{\LTpost}{3pt}
\setlength{\LTcapwidth}{\linewidth}

\begin{longtable}{@{}*{4}{
>{\centering\arraybackslash}p{\dimexpr 0.250000\linewidth-2\tabcolsep\relax}}@{}}
\caption{All integer solutions with negative balance: right triangles}
\label{tab:2}\\[-2pt]
\toprule
$2-\varepsilon$ & $t$ & $s$ & $r$ \\
\midrule
\endfirsthead
\caption[]{Continued}\\[-2pt]
\toprule
$2-\varepsilon$ & $t$ & $s$ & $r$ \\
\midrule
\endhead
\midrule
\endfoot
\bottomrule
\endlastfoot
$2$ & $0$ & $5$ & $0,1$ \\
$2$ & $1$ & $4$ & $0$ \\
\end{longtable}
\endgroup

Substituting the integer solutions in the table into the vertex angle equation gives exactly $\alpha\in\{18^\circ,(45/2)^\circ\}$. For every other parameter, all candidate vertex types satisfy $r-s\ge0$, so the zero-sum argument gives $r=s$ for every vertex type of positive limiting weight. The vertex angle equation then gives $s+t=2(2-\varepsilon)$, and hence $q\le8$. Substitution into inequality~(\ref{eq:degree-bound}) gives $N\le15$.

For family A, put $\mathcal E=3(2-\varepsilon)-s-2t$. The vertex angle equation is $\alpha(3r-2s-t)=180^\circ\mathcal E$. Thus $\mathcal E<0$ gives $s+2t-3(2-\varepsilon)>0$, and
\[
\alpha=\frac{180^\circ(s+2t-3(2-\varepsilon))}{2s+t-3r}<30^\circ
\quad\Longleftrightarrow\quad 3r+4s+11t<18(2-\varepsilon).
\]
Together with $s+2t>3(2-\varepsilon)$, this gives all the integer solutions in Table~\ref{tab:3}.

\begingroup
\tablefont
\setlength{\tabcolsep}{2pt}
\renewcommand{\arraystretch}{1}
\setlength{\LTpre}{2pt}
\setlength{\LTpost}{3pt}
\setlength{\LTcapwidth}{\linewidth}

\begin{longtable}{@{}*{8}{
>{\centering\arraybackslash}p{\dimexpr 0.125000\linewidth-2\tabcolsep\relax}}@{}}
\caption{All integer solutions with negative balance: family A}
\label{tab:3}\\[-2pt]
\toprule
$2-\varepsilon$ & $t$ & $s$ & $r$ & $2-\varepsilon$ & $t$ & $s$ & $r$ \\
\midrule
\endfirsthead
\caption[]{Continued}\\[-2pt]
\toprule
$2-\varepsilon$ & $t$ & $s$ & $r$ & $2-\varepsilon$ & $t$ & $s$ & $r$ \\
\midrule
\endhead
\midrule
\endfoot
\bottomrule
\endlastfoot
$1$ & $0$ & $4$ & $0$ & $2$ & $1$ & $5$ & $0,1$ \\
$2$ & $0$ & $7$ & $0,1,2$ & $2$ & $1$ & $6$ & $0$ \\
$2$ & $0$ & $8$ & $0,1$ & $2$ & $2$ & $3$ & $0$ \\
\end{longtable}
\endgroup

Substituting all integer solutions into the vertex angle equation gives the following set of exceptional values of $\alpha$:
\[
\{(90/7)^\circ,(180/11)^\circ,(45/2)^\circ,(360/13)^\circ\}.
\]
For every other parameter, every vertex type of positive limiting weight satisfies $\mathcal E=0$ and $3r-2s-t=0$. The vertex types satisfying these two equations are
\[
(4,6,0;0),(3,4,1;0),(2,2,2;0),
(1,0,3;0),(2,3,0;1),(1,1,1;1).
\]
Let the limiting weights of these six vertex types, in this order, be $x_1,x_2,x_3,x_4,x_5,x_6$. Since the incidence count for each of the three angles is $1$, we obtain $x_4=2x_1+x_2+x_5$ and $6x_1+4x_2+2x_3+3x_5+x_6=1$. Substitution into equation~(\ref{eq:angle-incidence}) gives
\[
N=16-8x_2-8x_3-16x_5-8x_6\le16.
\]

For family B, put $\mathcal E=2(2-\varepsilon)-s-t$. The vertex angle equation is $\alpha(2r-3s+t)=180^\circ\mathcal E$. Parameters admitting $\mathcal E<0$ must satisfy $s+t>2(2-\varepsilon)$, and
\[
\alpha=\frac{180^\circ(s+t-2(2-\varepsilon))}{3s-t-2r}<30^\circ
\quad\Longleftrightarrow\quad 2r+3s+7t<12(2-\varepsilon).
\]
All integer solutions of these conditions are listed in Table~\ref{tab:4}.

\begingroup
\tablefont
\setlength{\tabcolsep}{2pt}
\renewcommand{\arraystretch}{1}
\setlength{\LTpre}{2pt}
\setlength{\LTpost}{3pt}
\setlength{\LTcapwidth}{\linewidth}

\begin{longtable}{@{}*{8}{
>{\centering\arraybackslash}p{\dimexpr 0.125000\linewidth-2\tabcolsep\relax}}@{}}
\caption{All integer solutions with negative balance: family B}
\label{tab:4}\\[-2pt]
\toprule
$2-\varepsilon$ & $t$ & $s$ & $r$ & $2-\varepsilon$ & $t$ & $s$ & $r$ \\
\midrule
\endfirsthead
\caption[]{Continued}\\[-2pt]
\toprule
$2-\varepsilon$ & $t$ & $s$ & $r$ & $2-\varepsilon$ & $t$ & $s$ & $r$ \\
\midrule
\endhead
\midrule
\endfoot
\bottomrule
\endlastfoot
$1$ & $0$ & $3$ & $0,1$ & $2$ & $1$ & $4$ & $0,1,2$ \\
$2$ & $0$ & $5$ & $0,1,2,3,4$ & $2$ & $1$ & $5$ & $0$ \\
$2$ & $0$ & $6$ & $0,1,2$ & $2$ & $2$ & $3$ & $0$ \\
$2$ & $0$ & $7$ & $0,1$ & & & & \\
\end{longtable}
\endgroup

Substituting these integer solutions into the vertex angle equation gives exactly the following set of exceptional values of $\alpha$:
\[
\{12^\circ,(180/13)^\circ,(180/11)^\circ,20^\circ,(45/2)^\circ,(180/7)^\circ,(540/19)^\circ\}.
\]
For every other parameter, every vertex type of positive limiting weight satisfies $\mathcal E=0$ and must be one of the vertex types obtained from the six vertex types for family A above by interchanging $r,s$. The same argument gives $N\le16$.

For family C, put $\mathcal E=3(2-\varepsilon)-s-2t$. The vertex angle equation is $\alpha(r-s)=60^\circ\mathcal E$. Parameters admitting $\mathcal E<0$ must satisfy $s+2t>3(2-\varepsilon)$, and
\[
\alpha=\frac{60^\circ(s+2t-3(2-\varepsilon))}{s-r}<30^\circ
\quad\Longleftrightarrow\quad r+s+4t<6(2-\varepsilon).
\]
All integer solutions of these conditions are listed in Table~\ref{tab:5}.

\begingroup
\tablefont
\setlength{\tabcolsep}{2pt}
\renewcommand{\arraystretch}{1}
\setlength{\LTpre}{2pt}
\setlength{\LTpost}{3pt}
\setlength{\LTcapwidth}{\linewidth}

\begin{longtable}{@{}*{8}{
>{\centering\arraybackslash}p{\dimexpr 0.125000\linewidth-2\tabcolsep\relax}}@{}}
\caption{All integer solutions with negative balance: family C}
\label{tab:5}\\[-2pt]
\toprule
$2-\varepsilon$ & $t$ & $s$ & $r$ & $2-\varepsilon$ & $t$ & $s$ & $r$ \\
\midrule
\endfirsthead
\caption[]{Continued}\\[-2pt]
\toprule
$2-\varepsilon$ & $t$ & $s$ & $r$ & $2-\varepsilon$ & $t$ & $s$ & $r$ \\
\midrule
\endhead
\midrule
\endfoot
\bottomrule
\endlastfoot
$1$ & $0$ & $4$ & $0,1$ & $2$ & $0$ & $11$ & $0$ \\
$1$ & $0$ & $5$ & $0$ & $2$ & $1$ & $5$ & $0,1,2$ \\
$2$ & $0$ & $7$ & $0,1,2,3,4$ & $2$ & $1$ & $6$ & $0,1$ \\
$2$ & $0$ & $8$ & $0,1,2,3$ & $2$ & $1$ & $7$ & $0$ \\
$2$ & $0$ & $9$ & $0,1,2$ & $2$ & $2$ & $3$ & $0$ \\
$2$ & $0$ & $10$ & $0,1$ & & & & \\
\end{longtable}
\endgroup

Substituting the integer solutions in the table into the vertex angle equation gives exactly the following set of exceptional values of $\alpha$:
\[
\{(60/7)^\circ,10^\circ,12^\circ,15^\circ,(120/7)^\circ,20^\circ,(45/2)^\circ,24^\circ,(180/7)^\circ,(80/3)^\circ,(300/11)^\circ\}.
\]
For every other parameter, the zero-sum balance and the vertex angle equation force every vertex type of positive limiting weight to satisfy both $\mathcal E=0$ and $r=s$. The vertex types satisfying these two conditions are
\[
(6,6,0;0),(4,4,1;0),(2,2,2;0),
(0,0,3;0),(3,3,0;1),(1,1,1;1).
\]
\phantomsection\label{obs:c-family-bound}
Again denote their limiting weights, in this order, by $x_1,x_2,x_3,x_4,x_5,x_6$. Counting the $\alpha$- and $\gamma$-angles gives
\[
6x_1+4x_2+2x_3+3x_5+x_6=1,\qquad
x_2+2x_3+3x_4+x_6=1.
\]
Using these two equations to eliminate $x_1,x_4$ from equation~(\ref{eq:angle-incidence}) gives
\[
N=21-18x_2-18x_3-25x_5-13x_6\le21.
\tag{7}\label{eq:c-family-bound}
\]

For nonexceptional parameters, the balance at every actual vertex is a nonnegative integer. In family A, for example, each of the three angle-incidence counts in the disk is $f_R+O(R)$, and counting the total non-straight angle gives $\sum (2-\varepsilon)=f_R+O(R)$. Hence summing over vertices in $B_R$ gives
\[
\sum \mathcal E=3\sum (2-\varepsilon)-\sum s-2\sum t=O(R).
\]
Each vertex with $\mathcal E>0$ contributes at least $1$ to the left-hand side, so there are only $O(R)$ such vertices. The same counting argument applies to families B and C; for family R, use $\sum(r-s)=O(R)$. Since the vertex degrees are uniformly bounded, these vertices also contribute only $O(R)$ to the sum of the neighbor counts over tiles. Their limiting weights are therefore zero along every choice of radius subsequence.

We have now proved $N\le21$ for the nonexceptional parameters in all four families. The remaining parameters are determined by all the integer solutions in Tables~\ref{tab:2}--\ref{tab:5}, giving the $24$ exceptional angle triples listed in Table~\ref{tab:6}.

\begingroup
\tablefont
\setlength{\tabcolsep}{2pt}
\renewcommand{\arraystretch}{1}
\setlength{\LTpre}{2pt}
\setlength{\LTpost}{3pt}
\setlength{\LTcapwidth}{\linewidth}

\begin{longtable}{@{}*{8}{
>{\centering\arraybackslash}p{\dimexpr 0.125000\linewidth-2\tabcolsep\relax}}@{}}
\caption{The 24 exceptional angle triples (in degrees)}
\label{tab:6}\\[-2pt]
\toprule
Family & $\alpha$ & $\beta$ & $\gamma$ &
Family & $\alpha$ & $\beta$ & $\gamma$ \\
\midrule
\endfirsthead
\caption[]{Continued}\\[-2pt]
\toprule
Family & $\alpha$ & $\beta$ & $\gamma$ &
Family & $\alpha$ & $\beta$ & $\gamma$ \\
\midrule
\endhead
\midrule
\endfoot
\bottomrule
\endlastfoot
R & $18$ & $72$ & $90$ &
B & $540/19$ & $900/19$ & $1980/19$ \\
R & $45/2$ & $135/2$ & $90$ &
C & $60/7$ & $360/7$ & $120$ \\
A & $90/7$ & $360/7$ & $810/7$ &
C & $10$ & $50$ & $120$ \\
A & $180/11$ & $540/11$ & $1260/11$ &
C & $12$ & $48$ & $120$ \\
A & $45/2$ & $45$ & $225/2$ &
C & $15$ & $45$ & $120$ \\
A & $360/13$ & $540/13$ & $1440/13$ &
C & $120/7$ & $300/7$ & $120$ \\
B & $12$ & $72$ & $96$ &
C & $20$ & $40$ & $120$ \\
B & $180/13$ & $900/13$ & $1260/13$ &
C & $45/2$ & $75/2$ & $120$ \\
B & $180/11$ & $720/11$ & $1080/11$ &
C & $24$ & $36$ & $120$ \\
B & $20$ & $60$ & $100$ &
C & $180/7$ & $240/7$ & $120$ \\
B & $45/2$ & $225/4$ & $405/4$ &
C & $80/3$ & $100/3$ & $120$ \\
B & $180/7$ & $360/7$ & $720/7$ &
C & $300/11$ & $360/11$ & $120$ \\
\end{longtable}
\endgroup

All $24$ exceptional angle triples have rational angles, and none is the $30^\circ$-$60^\circ$-$90^\circ$ triangle. By the remark following Lemma~\ref{lem:unmatched}, no side of any of these prototiles is twice the length of another. The injection in Lemma~\ref{lem:unmatched} and the parity-correction inequality~(\ref{eq:parity-bound}) therefore apply.

\begin{lemma}[Weighted bounds and elimination of vertex types]
\label{lem:weighted-screening}
Let the prototile be a scalene triangle with no side twice the length of another. Take a finite list of candidate vertex types containing every vertex type that actually occurs. Then the following bounds and elimination rules hold.
\par\noindent\textup{(i)}~If rational coefficients $u,v,w,z$ satisfy
\[
q(q-2)\le ur+vs+wt+zP,\qquad z\ge0,
\tag{8}\label{eq:weighted-bound}
\]
for every vertex type in the list, then $N\le u+v+w$.
\par\noindent\textup{(ii)}~Let $K$ be a positive integer such that $N\le K$, and label the corner positions corresponding to $\alpha,\beta,\gamma$ by $1,2,3$, respectively. If an angle is absent from all vertex types in the list, no such tiling exists. Otherwise, let $\mu_i$ be the minimum of $q-2$ among the vertex types in the current list containing the $i$th angle. Every vertex type that actually occurs and contains the $i$th angle must satisfy
\[
q-2+\sum_{j\ne i}\mu_j\le K;
\]
if $\varepsilon=1$, it must also satisfy
\[
q-2+\mu_1+\mu_2+\mu_3\le K.
\]
Deleting candidate vertex types that fail these conditions and updating the $\mu_i$ at each round retains every vertex type that actually occurs, including those of zero limiting weight.
\end{lemma}

\begin{proof}
Summing inequality~(\ref{eq:weighted-bound}) with weights $x_\tau$ and using Lemma~\ref{lem:angle-incidence}\textup{(i)}, equation~(\ref{eq:angle-incidence}) and inequality~(\ref{eq:parity-bound}) gives
\[
N\le u+v+w+z\sum_\tau P(\tau)x_\tau
\le u+v+w.
\]
Since $N$ is an integer, the integer part of this upper bound may be used as $K$ in the next round of elimination.

To check inequality~(\ref{eq:weighted-bound}), choose a positive integer $L$ clearing the denominators of the coefficients, and write its slack as
\[
\Delta=L\bigl(ur+vs+wt+zP-q(q-2)\bigr).
\]
For fixed $\varepsilon,t,P$, the vertex angle equation gives
\[
r=\frac{180^\circ(2-\varepsilon)-s\beta-t\gamma}{\alpha}.
\]
Thus $r,q$ are affine functions of $s$, and the quadratic coefficient of $\Delta$, viewed as a polynomial in $s$, is $-L(\beta-\alpha)^2/\alpha^2\le0$. Hence $\Delta$ is concave. For each nonempty group, nonnegative slack at the smallest and largest feasible integer values of $s$ proves inequality~\eqref{eq:weighted-bound} for every vertex type in the group.

Fix an actual vertex and a tile contributing the $i$th angle there. This vertex contributes $q-2$ to equation~(1), while the other two corners of the tile contribute at least $\sum_{j\ne i}\mu_j$. Their sum therefore cannot exceed $K$. At a half-vertex, apply equation~(1) to the tile providing the straight angle there: its three corners contribute at least $\mu_1+\mu_2+\mu_3$, giving the second necessary condition. These conditions hold for each tile, so updating the $\mu_i$ on the candidate list and continuing the elimination retains all vertex types that actually occur, including those of zero limiting weight.
\end{proof}

Apart from $(18^\circ,72^\circ,90^\circ)$ and $(15^\circ,45^\circ,120^\circ)$, the remaining $22$ triples are treated by alternating the weighted bounds and elimination on individual tiles in Lemma~\ref{lem:weighted-screening}. The $31$ sets of coefficients are listed in Table~\ref{tab:7}. In each row, $(u,v,w,z)$ satisfies inequality~\eqref{eq:weighted-bound} for the candidate vertex types at that stage, and the last column gives the resulting upper bound $u+v+w$. For successive rows with the same angle triple, the preceding row's bound is used for elimination. Here $K$ is the positive integer upper bound used for elimination in that row, a dash means that no upper bound has yet been imposed, and the Types column gives the number of candidate vertex types at that stage. Family labels are as in Table~\ref{tab:1}.

\begingroup
\tablefont
\setlength{\tabcolsep}{2pt}
\renewcommand{\arraystretch}{1}
\setlength{\LTpre}{2pt}
\setlength{\LTpost}{3pt}
\setlength{\LTcapwidth}{\linewidth}

\begin{longtable}{@{}
>{\centering\arraybackslash}p{\dimexpr 0.055000\linewidth-2\tabcolsep\relax}
>{\centering\arraybackslash}p{\dimexpr 0.080000\linewidth-2\tabcolsep\relax}
>{\centering\arraybackslash}p{\dimexpr 0.035000\linewidth-2\tabcolsep\relax}
>{\centering\arraybackslash}p{\dimexpr 0.070000\linewidth-2\tabcolsep\relax}
>{\centering\arraybackslash}p{\dimexpr 0.145000\linewidth-2\tabcolsep\relax}
>{\centering\arraybackslash}p{\dimexpr 0.145000\linewidth-2\tabcolsep\relax}
>{\centering\arraybackslash}p{\dimexpr 0.145000\linewidth-2\tabcolsep\relax}
>{\centering\arraybackslash}p{\dimexpr 0.145000\linewidth-2\tabcolsep\relax}
>{\centering\arraybackslash}p{\dimexpr 0.180000\linewidth-2\tabcolsep\relax}@{}}
\caption{Exact weighted coefficients for the ordinary exceptional cases}
\label{tab:7}\\[-2pt]
\toprule
Fam. & $\alpha\ ({}^\circ)$ & $K$ & Types &
$u$ & $v$ & $w$ & $z$ & Upper bound for $N$ \\
\midrule
\endfirsthead
\caption[]{Continued}\\[-2pt]
\toprule
Fam. & $\alpha\ ({}^\circ)$ & $K$ & Types &
$u$ & $v$ & $w$ & $z$ & Upper bound for $N$ \\
\midrule
\endhead
\midrule
\endfoot
\bottomrule
\endlastfoot
R & $45/2$ & \textendash{} & $23$ &
$14$ & $3$ & $2$ & $1$ & $19$ \\
A & $90/7$ & \textendash{} & $23$ &
$26$ & $104/31$ & $-305/31$ & $357/31$ & $605/31$ \\
A & $180/11$ & \textendash{} & $24$ &
$20$ & $48/11$ & $-78/11$ & $102/11$ & $190/11$ \\
A & $45/2$ & \textendash{} & $27$ &
$14$ & $6$ & $-9/2$ & $15/2$ & $31/2$ \\
A & $360/13$ & \textendash{} & $15$ &
$140/13$ & $80/13$ & $-38/13$ & $6$ & $14$ \\
B & $12$ & \textendash{} & $20$ &
$28$ & $-4/11$ & $-31/11$ & $185/11$ & $273/11$ \\
B & $12$ & $24$ & $17$ &
$3884/195$ & $48/65$ & $-263/195$ & $147/13$ & $251/13$ \\
B & $180/13$ & \textendash{} & $20$ &
$24$ & $-28/11$ & $-8/11$ & $140/11$ & $228/11$ \\
B & $180/11$ & \textendash{} & $21$ &
$20$ & $5/4$ & $-13/8$ & $93/8$ & $157/8$ \\
B & $20$ & \textendash{} & $23$ &
$16$ & $32/5$ & $-16/5$ & $56/5$ & $96/5$ \\
B & $45/2$ & \textendash{} & $12$ &
$14$ & $29/13$ & $-8/13$ & $99/13$ & $203/13$ \\
B & $180/7$ & \textendash{} & $26$ &
$12$ & $4$ & $-1$ & $7$ & $15$ \\
B & $540/19$ & \textendash{} & $9$ &
$188/19$ & $98/19$ & $-20/19$ & $6$ & $14$ \\
C & $60/7$ & \textendash{} & $23$ &
$40$ & $41/31$ & $-235/31$ & $798/31$ & $1046/31$ \\
C & $60/7$ & $33$ & $21$ &
$14719/462$ & $325/154$ & $-379/66$ & $445/22$ & $621/22$ \\
C & $60/7$ & $28$ & $19$ &
$569/21$ & $18/7$ & $-14/3$ & $17$ & $25$ \\
C & $60/7$ & $25$ & $17$ &
$2166/91$ & $263/91$ & $-51/13$ & $192/13$ & $296/13$ \\
C & $60/7$ & $22$ & $15$ &
$173/9$ & $10/3$ & $-26/9$ & $35/3$ & $59/3$ \\
C & $10$ & \textendash{} & $23$ &
$34$ & $-25/31$ & $-172/31$ & $609/31$ & $857/31$ \\
C & $10$ & $27$ & $20$ &
$437/18$ & $25/18$ & $-11/3$ & $14$ & $22$ \\
C & $10$ & $22$ & $17$ &
$859/39$ & $74/39$ & $-42/13$ & $165/13$ & $269/13$ \\
C & $12$ & \textendash{} & $24$ &
$28$ & $25/11$ & $-53/11$ & $192/11$ & $280/11$ \\
C & $12$ & $25$ & $22$ &
$2851/120$ & $319/120$ & $-185/48$ & $233/16$ & $361/16$ \\
C & $12$ & $22$ & $20$ &
$944/45$ & $131/45$ & $-29/9$ & $38/3$ & $62/3$ \\
C & $120/7$ & \textendash{} & $14$ &
$166/9$ & $37/9$ & $-26/9$ & $35/3$ & $59/3$ \\
C & $20$ & \textendash{} & $29$ &
$16$ & $6$ & $-11/4$ & $45/4$ & $77/4$ \\
C & $45/2$ & \textendash{} & $11$ &
$14$ & $81/13$ & $-30/13$ & $129/13$ & $233/13$ \\
C & $24$ & \textendash{} & $17$ &
$38/3$ & $8$ & $-29/12$ & $41/4$ & $73/4$ \\
C & $180/7$ & \textendash{} & $12$ &
$12$ & $8$ & $-9/4$ & $39/4$ & $71/4$ \\
C & $80/3$ & \textendash{} & $9$ &
$373/36$ & $347/36$ & $-9/4$ & $39/4$ & $71/4$ \\
C & $300/11$ & \textendash{} & $8$ &
$453/44$ & $427/44$ & $-9/4$ & $39/4$ & $71/4$ \\
\end{longtable}
\endgroup

Every branch in Table~\ref{tab:7} ends with $N<22$, so these $22$ exceptional triples all satisfy $N\le21$. The $581$ inequalities on individual vertex types reduce, by the preceding concavity argument for the slacks, to $426$ endpoint substitutions; the endpoint slacks and successive deletion records are given in Appendix Tables~\ref{tab:8}--\ref{tab:38}. The remaining two triples require additional local geometric counts.

\par\addvspace{6.5bp}\noindent
\begin{minipage}{\linewidth}
\begin{proposition}\label{prop:18}
For the $18^\circ$-$72^\circ$-$90^\circ$ triangular prototile, the common neighbor number cannot satisfy $N\ge22$.
\end{proposition}
\end{minipage}\par

\begin{proof}
The vertex angle equation
\[
r+4s+5t=20-10\varepsilon
\]
has $22$ candidate vertex types. Put $\tau_B=(0,5,0;0)$ and $b_0=x_{\tau_B}$. Inequality~\eqref{eq:parity-bound} gives $b_0\le h$. Write $\mathbf1_{\tau_B}$ for the indicator of type $\tau_B$. Every candidate vertex type satisfies
\[
q(q-2)\le18r+\frac{14}{5}s+2t-\varepsilon+\mathbf1_{\tau_B}.
\]
The corresponding integer slacks are listed in the $\Delta_0$ column of Appendix Table~\ref{tab:39}. Summing with the limiting weights and using Lemma~\ref{lem:angle-incidence} and the preceding bound for $b_0$ gives
\[
N\le\frac{114}{5}-h+b_0\le\frac{114}{5}<23.
\]
Since $N$ is an integer, $N\le22$. Suppose that $N=22$.

Put $\tau_C=(0,0,4;0)$ and $\tau_R=(0,0,2;1)$, with $c_0=x_{\tau_C}$ and $r_0=x_{\tau_R}$. Let the number of $72^\circ$ angles at vertices not of type $\tau_B$ be $E$, the number of right angles at vertices of neither type $\tau_C$ nor type $\tau_R$ be $M$, and the number of straight angles at vertices not of type $\tau_R$ be $U$. All these are limits of counts normalized by the number of tiles and taken along the fixed subsequence. Counting the $72^\circ$ angles and the right angles separately gives
\[
5b_0+E=1,\qquad 4c_0+2r_0+M=1.
\tag{9}\label{eq:special-angle-counts}
\]

For a tile $T$, let $n_R(T)$ and $n_U(T)$ be its numbers of subdivision vertices of type $\tau_R$ and of types other than $\tau_R$, respectively. Two consecutive subdivision vertices on a side cannot both have type $\tau_R$, since the tile on the opposite side of the segment between them would have two right angles. Hence a side with $u$ subdivision vertices not of type $\tau_R$ has at most $u+1$ subdivision vertices of type $\tau_R$. If the tile across one of $T$'s legs also provides a right angle at $T$'s right-angled endpoint, the first subdivision vertex cannot have type $\tau_R$, and the bound decreases to $u$. When the right angle of $T$ is at a vertex of type $\tau_C$, this stronger restriction applies to both legs. Summing over its three sides gives
\[
n_R(T)\le n_U(T)+1.
\]
When the right angle of $T$ is at a vertex of type $\tau_R$, the bound $u$ applies to one leg, giving
\[
n_R(T)\le n_U(T)+2.
\]

Now let $T$ be a tile whose $\beta$-corner is at a vertex of type $\tau_B$. Let $q_\alpha,q_\beta,q_\gamma$ be the degrees at its three corners, and let $\mathcal V_\pi(T)$ be its set of subdivision vertices. Substituting $q_\beta=5$ and $N=22$ into equation~(1) gives
\[
q_\alpha+q_\gamma+
\sum_{v\in\mathcal V_\pi(T)}(q(v)-2)=23.
\tag{10}\label{eq:subdivided-tile-neighbors}
\]
For $\varepsilon=0$, the vertex angle equation gives $q=20-3s-4t$, so $q$ cannot be $18$ or $19$; for $\varepsilon=1$, we have $q\le11$. Thus a vertex containing an $\alpha=18^\circ$ angle has degree either $20$ or at most $17$. Suppose that this tile satisfies $n_U(T)=0$. If its $\gamma$-corner is at a vertex of type $\tau_C$, the total contribution of its subdivision vertices to equation~\eqref{eq:subdivided-tile-neighbors} is $n_R(T)\le1$, forcing $q_\alpha$ to be $18$ or $19$, contrary to the degree restriction from the vertex angle equation. If its $\gamma$-corner is at a vertex of type $\tau_R$, then $n_R(T)\le2$, and equation~\eqref{eq:subdivided-tile-neighbors} can hold only when $q_\alpha=20$ and $n_R(T)=0$. In this case, $T$ has no subdivision vertices.

At each vertex of type $\tau_B$, choose a ray with unequal side lengths as in Lemma~\ref{lem:unmatched}, and select the tile with the longer side. This tile has its $\beta$-corner at that vertex, and its longer side has a subdivision vertex. Since each tile has only one $\beta$-corner, the selected tiles are distinct. The local restrictions just proved imply that each selected tile has its $\gamma$-corner counted by $M$ or has at least one subdivision vertex counted by $U$. Each point counted by $U$ has a unique straight-angle tile, so this correspondence is injective. The disk-counting framework of Section~\ref{sec:basic-facts} gives
\[
b_0\le M+U.
\]
Similarly, every tile whose $\gamma$-corner is at a vertex of type $\tau_C$ either has its $\beta$-corner outside type $\tau_B$ or has at least one subdivision vertex counted by $U$. Hence
\[
4c_0\le E+U.
\]
Equations~\eqref{eq:special-angle-counts} now give
\[
b_0\le M+U,\qquad 5b_0\le2r_0+M+U.
\tag{11}\label{eq:special-geometric-counts}
\]

To apply the counting constraints~\eqref{eq:special-geometric-counts}, define the deficit of each vertex type by
\[
\delta(\tau)=18r+3s+2t-q(q-2).
\]
Its sum with the limiting weights is $23-N=1$. The pure $\alpha$ vertex type of degree $20$ and vertex types $\tau_B$ and $\tau_C$ have deficit zero, while vertex type $\tau_R$ has deficit $1$. For fixed $\varepsilon,t$, the slack $8\delta(\tau)-9s-29t-29\varepsilon$ is a concave quadratic in $s$. Checking the feasible endpoints as in the proof of Lemma~\ref{lem:weighted-screening} gives nonnegative slacks for all candidate vertex types other than $(20,0,0;0)$, $\tau_B$, $\tau_C$ and $\tau_R$; the exact values are listed in the $\Delta_1$ column of Appendix Table~\ref{tab:39}. Thus these vertex types satisfy
\[
8\delta(\tau)\ge9s+29t+29\varepsilon.
\tag{12}\label{eq:deficit-estimate}
\]
Outside the four vertex types above, the weighted sums of $s,t,\varepsilon$ are $E,M,U$, respectively. Summing inequality~\eqref{eq:deficit-estimate} with the limiting weights and adding the contribution $r_0$ from type $\tau_R$ gives
\[
1\ge r_0+\frac98E+\frac{29}{8}(M+U).
\]
By inequalities~\eqref{eq:special-geometric-counts},
\[
r_0+\frac{29}{8}(M+U)-\frac{45}{8}b_0
=\frac{25}{8}(M+U-b_0)
+\frac12(2r_0+M+U-5b_0)\ge0.
\]
\begin{samepage}
Thus
\[
1\ge\frac98(E+5b_0)=\frac98,
\]
a contradiction.
\end{samepage}
\end{proof}

\par\addvspace{6.5bp}\noindent
\begin{minipage}{\linewidth}
\begin{proposition}\label{prop:15}
For the $15^\circ$-$45^\circ$-$120^\circ$ triangular prototile, $N\le21$.
\end{proposition}
\end{minipage}\par

\begin{proof}
Suppose that $N\ge22$. Write $\alpha=15^\circ$, $\beta=45^\circ$, $\gamma=120^\circ$, and let $a<b<c$ be their respective opposite side lengths. The vertex angle equation
\[
r+3s+8t=24-12\varepsilon
\]
gives $26$ candidate vertex types, listed in Case C24 of Appendix Table~\ref{app:rigid-slacks-combined}. The only degree-$3$ type is $(0,0,3;0)$, and the only types of degrees $4$ and $5$ are $\tau_X=(1,1,1;1)$ and $\tau_Y=(0,4,0;1)$, respectively. The only types of degree greater than $18$ are $(24,0,0;0)$, $\tau_Z=(21,1,0;0)$ and $\tau_W=(18,2,0;0)$, of degrees $24,22,20$, respectively. In particular, degrees $21$ and $23$ do not occur.

For any tile $T$, write $q_\alpha,q_\beta,q_\gamma$ for its three corner degrees. Suppose that $T$ has $k\ge1$ subdivision vertices. Its $\alpha$- and $\beta$-corners and each subdivision vertex have degree at least $4$, while $q_\gamma\ge3$. If $k=1$, equation~(1) gives a total degree of $N+8$ at the four boundary vertices. Applying the reciprocal-sum estimate in the proof of Lemma~\ref{lem:local-curvature}, we obtain
\[
\kappa(T)\le-1+\frac13+\frac24+\frac1{(N+8)-3-4-4}
=-\frac16+\frac1{N-3}.
\]
If $k\ge2$, the same lower bounds on the degrees give $\kappa(T)\le1/3-k/4\le-1/6$. Thus every tile with subdivision vertices satisfies
\[
\kappa(T)\le-\frac16+\frac1{N-3}<0.
\]
Consequently, every tile of nonnegative curvature has no subdivision vertices. For a tile without subdivision vertices, equation~(1) and the curvature definition give, respectively,
\[
q_\alpha+q_\beta+q_\gamma=N+6,\qquad
\kappa(T)=-\frac12+\frac1{q_\alpha}+\frac1{q_\beta}+\frac1{q_\gamma}.
\tag{13}\label{eq:triangle-tile-curvature}
\]

At each degree-$3$ vertex, all three incident tiles contribute a $\gamma$-angle. By Lemma~\ref{lem:unmatched}, at least one of these tiles has a subdivision vertex and hence negative curvature. Fix one such tile at each of these vertices, and transfer to it the curvature of each other incident tile of nonnegative curvature. Since a receiving tile has only one $\gamma$-corner, it receives at most two transfers in the whole tiling.

\emph{(i) $N\ge26$.} Let $T$ have nonnegative curvature, so it has no subdivision vertices. By equation~\eqref{eq:triangle-tile-curvature}, the sum of the reciprocals of its three corner degrees is at least $1/2$. If $q_\gamma\ge4$, then $4\le q_\alpha,q_\beta,q_\gamma\le24$, and equation~\eqref{eq:triangle-tile-curvature} gives $q_\alpha+q_\beta+q_\gamma=N+6\ge32$. If two degrees are $4$, the third is at least $24$. Hence this occurs only when $N=26$ and $(q_\alpha,q_\beta,q_\gamma)=(24,4,4)$. In every other case, the two smallest entries are at least $4,5$, respectively. By the argument in the proof of Lemma~\ref{lem:local-curvature}, reduce them to $4,5$ and add the removed amounts to the third entry, which then becomes at least $32-4-5=23$. Hence
\[
\frac1{q_\alpha}+\frac1{q_\beta}+\frac1{q_\gamma}
\le\frac14+\frac15+\frac1{23}<\frac12,
\]
so $\kappa(T)<0$, a contradiction.

If $(q_\alpha,q_\beta,q_\gamma)=(24,4,4)$, let $v,u$ be the $\alpha$- and $\beta$-corners of $T$, respectively, so that $vu$ is its longest side. Since a degree-$24$ vertex contains only $\alpha$-angles, $v$ is a full vertex, and the neighboring tile $T_1$ across $vu$ also has its $\alpha$-corner at $v$. As $T$ has no subdivision vertices and $c$ is the longest side, $T_1$ shares the entire side $vu$ with $T$. The endpoints of a longest side are the $\alpha$- and $\beta$-corners, so $T_1$ contributes a $\beta$-angle at $u$. The degree-$4$ vertex $u$ would then contain two $\beta$-angles, contrary to type $\tau_X$. Therefore every nonnegative-curvature tile has $q_\gamma=3$.

Substituting $q_\gamma=3$ into equation~\eqref{eq:triangle-tile-curvature} gives
\[
\kappa(T)=-\frac16+\frac1{q_\alpha}+\frac1{q_\beta}.
\]
Since $q_\alpha+q_\beta\ge29$ and $q_\alpha,q_\beta\le24$, both degrees are at least $5$. If one degree is $5$, then, as degree-$5$ vertices contain only $\beta$-angles, we have $N=26$ and $(q_\alpha,q_\beta,q_\gamma)=(24,5,3)$, with curvature $3/40$. For any other tile, if the smaller degree is $6$, the other is at least $23$. Since degree $23$ does not occur, the other degree must be $24$, and equation~\eqref{eq:triangle-tile-curvature} gives $\kappa(T)\le1/24$. If both degrees are at least $7$, the preceding reciprocal-sum estimate and equation~\eqref{eq:triangle-tile-curvature} give
\[
\kappa(T)\le-\frac16+\frac17+\frac1{22}<\frac1{24}.
\]
Thus every nonnegative-curvature tile other than type $(24,5,3)$ has curvature at most $1/24$.

For a tile $T$ with corner degrees $(24,5,3)$, let $T_1$ be the neighboring tile across its $b$ side. Since both endpoints are full vertices and the side has no subdivision vertices, the two tiles share the entire $b$ side. The endpoints of degrees $24,3$ are the $\alpha$- and $\gamma$-corners of $T_1$, respectively. Equation~(1) gives
\[
q_\beta(T_1)+\sum_{x\in\mathcal V_\pi(T_1)}(q(x)-2)=5.
\]
The first term is at least $4$, while the sum is zero or at least $2$. Therefore $T_1$ has no subdivision vertices and $q_\beta(T_1)=5$, so it also has corner degrees $(24,5,3)$. Let $u,u'$ be the $\beta$-corners of $T,T_1$, respectively, let $w$ be their common $\gamma$-corner, and let $T_2$ be the third tile at $w$. Since $T_2$ contributes a $\gamma$-angle at $w$, interchanging $T,T_1$ if necessary gives the local configuration in Fig.~\ref{fig:curvature-local}(a). Thus $q_\beta(T_2)=5$, and $u'$ is a degree-$5$ subdivision vertex of $T_2$. If this were its only subdivision vertex, equation~(1) would give $q_\alpha(T_2)=21$, contrary to the preceding degree restriction. Hence $T_2$ has at least two subdivision vertices and curvature at most $-1/6$.

\begin{samepage}
Since $T,T_1$ both have positive curvature, the receiving tile at $w$ must be $T_2$, and it receives curvature only from these two tiles. Every other receiving tile receives at most two contributions, each at most $1/24$. After the transfers, the curvature of $T_2$ and that of any other receiving tile are bounded above, respectively, by
\[
-\frac16+2\cdot\frac3{40}=-\frac1{60}<0,
\qquad
-\frac16+\frac1{23}+2\cdot\frac1{24}<0.
\]\par
\end{samepage}

\begin{samepage}
All initially nonnegative-curvature tiles now have curvature zero, and negative-curvature tiles receiving no transfer are unchanged. For fixed $N$, the original negative curvatures take only finitely many values, since $q\le24$ and $m\le N$. Thus all initially negative-curvature tiles retain a uniform negative upper bound after the transfers. Let $f_R^-$ count the initially negative-curvature tiles in $\mathcal F_R$. Every initially nonnegative-curvature tile has been assigned to a negative-curvature tile, and each negative-curvature tile is assigned at most two such tiles. By the disk count in Sect.~\ref{sec:basic-facts},
\[
f_R\le3f_R^-+O(R).
\]
Together with the tile-count estimate in Sect.~\ref{sec:basic-facts}, this gives a uniformly positive proportion of negative-curvature tiles. Consequently, the average curvature after the transfers has a uniform negative upper bound, contradicting the average-curvature conclusion of Lemma~\ref{lem:curvature}. Hence $N\le25$.\par
\end{samepage}

\begin{figure}[htbp]
\centering
\includegraphics[width=110mm]{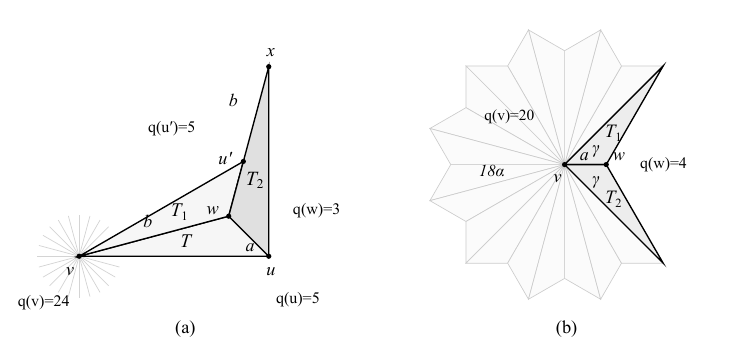}
\caption{Local configurations}\label{fig:curvature-local}
\end{figure}

Suppose that a pure $\alpha$ vertex of degree $24$ occurs. An incident tile with a subdivision vertex would give $N\ge27$ by equation~(1), so all incident tiles have no subdivision vertices. The radial $c$ sides therefore match completely. Their common outer endpoints contain two $\beta$-angles and hence have degree at least $5$, so equation~(1) gives $N\ge26$, contradicting $N\le25$.

If a vertex of type $\tau_Z$ occurs, let $T_1$ be the unique tile contributing a $\beta$-angle there. Across its radial $a$ side lies a $b$ or $c$ side of a tile contributing an $\alpha$-angle at the central vertex. Hence the $\gamma$-corner of $T_1$ lies in the interior of a side of that neighboring tile. The neighboring tile has a subdivision vertex, and equation~(1) gives at least $25$ neighbors, so only $N=25$ remains to be considered. On the other hand, the $\gamma$-corner of $T_1$ is a half-vertex, so the angle equation allows only degrees $4$ and $6$ there. Also, $q_\alpha(T_1)\ge4$ and $q_\alpha(T_1)\ne5$, and each subdivision vertex contributes at least $2$. Substituting the central degree $22$ into equation~(1), we obtain
\[
N=16+q_\alpha(T_1)+q_\gamma(T_1)
+\sum_{x\in\mathcal V_\pi(T_1)}(q(x)-2).
\]
The right-hand side has minimum $24$, and increasing any term above its minimum increases the sum by at least $2$. It therefore cannot equal $25$. With these two vertex types excluded, every vertex has degree at most $20$.

\emph{(ii) $23\le N\le25$.} For a tile without subdivision vertices and with $q_\gamma\ge4$, we have $4\le q_\alpha,q_\beta,q_\gamma\le20$, and equation~\eqref{eq:triangle-tile-curvature} gives $q_\alpha+q_\beta+q_\gamma\ge29$. If two degrees were $4$, the third would be at least $21$, exceeding this upper bound. In the remaining cases, the reciprocal-sum estimate used in case~(i) gives
\[
\frac1{q_\alpha}+\frac1{q_\beta}+\frac1{q_\gamma}
\le\frac14+\frac15+\frac1{20}=\frac12.
\]
By equation~\eqref{eq:triangle-tile-curvature}, nonnegative curvature requires equality. Equality holds only when $N=23$ and $(q_\alpha,q_\beta,q_\gamma)=(20,5,4)$. Every other nonnegative-curvature tile has $q_\gamma=3$. In this case $q_\alpha+q_\beta\ge26$, and both degrees are at most $20$, so both are at least $6$. Its curvature is therefore at most $-1/6+1/6+1/20=1/20$. After the transfers at degree-$3$ vertices, every receiving tile has curvature at most
\[
-\frac16+\frac1{20}+2\cdot\frac1{20}=-\frac1{60}<0.
\]

Let $T$ be a zero-curvature tile with corner degrees $(20,5,4)$, and let $v,w$ be its $\alpha$- and $\gamma$-corners, respectively. Then $vw$ has length $b$, and $v$ is a full vertex of type $\tau_W$. The neighboring tile across $vw$ therefore also has a corner at $v$. If its side along $vw$ had length $a$, then $T$ would have a subdivision vertex, contradicting its zero curvature. If that side had length $b$, the neighboring tile would contribute an $\alpha$-angle at $v$, since no $\gamma$-angle occurs there. It would then contribute a $\gamma$-angle at $w$, contrary to the degree-$4$ type $\tau_X$. Thus the side has length $c$, and $w$ is a subdivision vertex of the neighboring tile, which has negative curvature. Assign $T$ to this negative-curvature tile, transferring zero curvature. Each negative-curvature tile is adjacent to at most $N$ such zero-curvature tiles and is assigned at most two further nonnegative-curvature tiles at degree-$3$ vertices. Hence
\[
f_R\le(N+3)f_R^-+O(R).
\]
Negative-curvature tiles still occupy a uniformly positive proportion, and the negative upper bound after the transfers is unchanged. The curvature-sum argument in case~(i) now contradicts Lemma~\ref{lem:curvature}. Therefore $N\le22$.

\emph{(iii) $N=22$.} Consider a tiling with $N\le22$ and a vertex $v$ of type $\tau_W$. If any incident tile had a subdivision vertex, equation~(1) would give $N\ge23$. Thus all $20$ incident tiles have no subdivision vertices, and all radial sides match completely. The two tiles $T_1,T_2$ contributing $\beta$-angles provide the only radial $a$ sides, so they share an entire $a$ side. Let $w$ be its other endpoint. Both tiles contribute a $\gamma$-angle at $w$. For either tile, equation~(1) with $q_\beta=20$ gives $q_\alpha+q(w)=N-14\le8$. Since $q_\alpha\ge4$, $q_\alpha\ne5$ and $q(w)\ge3$, we have $q_\alpha=4$ and $q(w)\in\{3,4\}$. If $q(w)=4$, Fig.~\ref{fig:curvature-local}(b) forces two $\gamma$-angles at $w$, contrary to the angle multiplicities of type $\tau_X$. If $q(w)=3$, the third tile at $w$ also contributes a $\gamma$-angle and has incident sides of lengths $a,b$, whereas the other sides of the first two tiles issuing from $w$ both have length $b$. Since $a<b$, the endpoint of the $a$ side of the third tile lies in the interior of a $b$ side of one of the first two tiles, contradicting the absence of subdivision vertices there. Hence no vertex of type $\tau_W$ occurs when $N\le22$. Taking $N=22$, every vertex has degree at most $18$.

For a tile without subdivision vertices and with $q_\gamma\ge4$, we have $4\le q_\alpha,q_\beta,q_\gamma\le18$, and equation~\eqref{eq:triangle-tile-curvature} gives $q_\alpha+q_\beta+q_\gamma=28$. Two degrees cannot both be $4$, nor can they be $4$ and $5$, since the third would exceed $18$. The same pairwise argument gives the list $(4,6,18)$ when the smallest degree is $4$, and $(5,5,18)$ otherwise. The latter gives a smaller reciprocal sum, so
\[
\frac1{q_\alpha}+\frac1{q_\beta}+\frac1{q_\gamma}
\le\frac14+\frac16+\frac1{18}<\frac12.
\]
Equation~\eqref{eq:triangle-tile-curvature} gives $\kappa(T)<0$. Hence every nonnegative-curvature tile has $q_\gamma=3$ and $q_\alpha+q_\beta=25$. Since both degrees are at most $18$, both are at least $7$, giving curvature at most
\[
\kappa(T)\le-\frac16+\frac17+\frac1{18}=\frac2{63}.
\]
Applying the transfer rule at degree-$3$ vertices once more, every receiving tile has curvature at most
\[
-\frac16+\frac1{19}+2\cdot\frac2{63}<0.
\]
The disk count and curvature-sum argument in case~(i) give a contradiction. Thus $N=22$ is impossible, and $N\le21$.
\end{proof}

The nonexceptional parameters in all four families, the $22$ ordinary exceptional triples and the two special triples all satisfy $N\le21$. Together with the case $\alpha\ge30^\circ$ and Proposition~\ref{prop:polygons}, this proves the upper-bound part of Theorem~\ref{thm:main}.

\subsection{The equality case \texorpdfstring{$N=21$}{N=21}}
\label{sec:rigidity}

\begin{proof}[Equality case of Theorem~\ref{thm:main}]
Fig.~\ref{fig:construction} gives a construction with $N=21$. Suppose that $N=21$. Proposition~\ref{prop:polygons} implies that the prototile is a triangle. If $\alpha<30^\circ$, the four-family reduction in Section~\ref{sec:triangles} applies. The bounds $N\le15$ for nonexceptional parameters in family R and $N\le16$ for nonexceptional parameters in families A and B leave only the nonexceptional parameters in family C and the $24$ exceptional cases in Table~\ref{tab:6}.

For the nonexceptional parameters in family C, retain the limiting weights $x_1,\ldots,x_6$ of the six zero-balance vertex types from Section~\ref{sec:triangles}. Equation~\eqref{eq:c-family-bound}, the nonnegativity of the weights, and $N=21$ give $x_2=x_3=x_5=x_6=0$. The only possible positive-weight vertex types are $(6,6,0;0)$ and $(0,0,3;0)$. The $\beta$- and $\gamma$-angle counts in Lemma~\ref{lem:angle-incidence}\textup{(i)} give $6x_1=1$ and $3x_4=1$, respectively, so $x_1=1/6$ and $x_4=1/3$. Both types are full vertices, so $h=0$. The angle multiplicities of the type corresponding to $x_4$ do not all have the same parity. Lemma~\ref{lem:unmatched} bounds the total weight of such vertices by $2h=0$, contradicting $x_4=1/3$. 

For the $24$ exceptional cases in Table~\ref{tab:6}, Section~\ref{sec:triangles} has ruled out a side twice the length of another, so inequality~\eqref{eq:parity-bound} applies. For the $22$ ordinary exceptions, apply Lemma~\ref{lem:weighted-screening}\textup{(ii)} with $K=21$, starting from all candidate vertex types and updating $\mu_i$ after each elimination round. Each final list satisfies inequality~\eqref{eq:weighted-bound} with the coefficients in Table~\ref{article-rigid-table-1}. The upper bounds in the last column are all strictly less than $21$, excluding these $22$ prototiles. The $378$ slacks and the vertex types deleted in each round are given in Cases C1--C22 of Appendix Table~\ref{app:rigid-slacks-combined}.

\begingroup
\tablefont
\setlength{\tabcolsep}{3pt}
\renewcommand{\arraystretch}{1}
\setlength{\LTcapwidth}{\textwidth}

\begin{longtable}{*{7}{
>{\centering\arraybackslash}p{\dimexpr 0.1428\textwidth-2\tabcolsep\relax}}}
\caption{Weighted coefficients for the ordinary exceptions in the $21$-neighbor case}
\label{article-rigid-table-1}\\
\toprule
Family & $\alpha\ ({}^\circ)$ & $u$ & $v$ & $w$ & $z$ & $u+v+w$ \\
\midrule
\endfirsthead
\caption[]{\textit{Continued}}\\
\toprule
Family & $\alpha\ ({}^\circ)$ & $u$ & $v$ & $w$ & $z$ & $u+v+w$ \\
\midrule
\endhead
\midrule
\multicolumn{7}{r}{\textit{Continued on the next page}}\\
\endfoot
\bottomrule
\endlastfoot
R & $45/2$ & $14$ & $3$ & $2$ & $1$ & $19$ \\
A & $90/7$ & $26$ & $104/31$ & $-305/31$ &
$357/31$ & $605/31$ \\
A & $180/11$ & $20$ & $48/11$ & $-78/11$ &
$102/11$ & $190/11$ \\
A & $45/2$ & $14$ & $6$ & $-9/2$ &
$15/2$ & $31/2$ \\
A & $360/13$ & $140/13$ & $80/13$ & $-38/13$ &
$6$ & $14$ \\
B & $12$ & $3884/195$ & $48/65$ & $-263/195$ &
$147/13$ & $251/13$ \\
B & $180/13$ & $24$ & $-28/11$ & $-8/11$ &
$140/11$ & $228/11$ \\
B & $180/11$ & $20$ & $5/4$ & $-13/8$ &
$93/8$ & $157/8$ \\
B & $20$ & $16$ & $32/5$ & $-16/5$ &
$56/5$ & $96/5$ \\
B & $45/2$ & $14$ & $29/13$ & $-8/13$ &
$99/13$ & $203/13$ \\
B & $180/7$ & $12$ & $4$ & $-1$ &
$7$ & $15$ \\
B & $540/19$ & $188/19$ & $98/19$ & $-20/19$ &
$6$ & $14$ \\
C & $60/7$ & $173/9$ & $10/3$ & $-26/9$ &
$35/3$ & $59/3$ \\
C & $10$ & $859/39$ & $74/39$ & $-42/13$ &
$165/13$ & $269/13$ \\
C & $12$ & $944/45$ & $131/45$ & $-29/9$ &
$38/3$ & $62/3$ \\
C & $120/7$ & $166/9$ & $37/9$ & $-26/9$ &
$35/3$ & $59/3$ \\
C & $20$ & $16$ & $6$ & $-11/4$ &
$45/4$ & $77/4$ \\
C & $45/2$ & $14$ & $81/13$ & $-30/13$ &
$129/13$ & $233/13$ \\
C & $24$ & $38/3$ & $8$ & $-29/12$ &
$41/4$ & $73/4$ \\
C & $180/7$ & $12$ & $8$ & $-9/4$ &
$39/4$ & $71/4$ \\
C & $80/3$ & $373/36$ & $347/36$ & $-9/4$ &
$39/4$ & $71/4$ \\
C & $300/11$ & $453/44$ & $427/44$ & $-9/4$ &
$39/4$ & $71/4$ \\
\end{longtable}
\endgroup

\begin{samepage}
For the $18^\circ$-$72^\circ$-$90^\circ$ triangle, let $v$ be a full vertex. Assume that every tile contributing an $\alpha$-angle at $v$ has no subdivision vertices and has its $\beta$-corner at a vertex of type $\tau_X=(1,1,1;1)$. Choose such a tile $T$, and let $u$ be its $\beta$-corner, so that $vu$ is its longest side. Let $T_1$ be the neighboring tile across $vu$. Since $v$ is a full vertex, the side $vu$ of $T$ has no subdivision vertices, and the prototile has no longer side, $T_1$ and $T$ share the entire side $vu$. If $T_1$ contributed an $\alpha$-angle at $v$, Fig.~\ref{fig:longest-pairing}(a) would give a $\beta$-angle from each of $T,T_1$ at $u$, contrary to the angle multiplicities of type $\tau_X$. Hence $T_1$ contributes a $\beta$-angle at $v$, giving the pairing in Fig.~\ref{fig:longest-pairing}(b). Each tile $T_1$ has only one longest side, and there is only one tile across that entire side. The pairing is therefore injective, so the multiplicities of the $\alpha$- and $\beta$-angles satisfy $r\le s$.\par
\end{samepage}

\par\addvspace{\intextsep}\noindent\begin{minipage}{\linewidth}
\centering
\makeatletter\def\@captype{figure}\makeatother
\includegraphics[width=115mm]{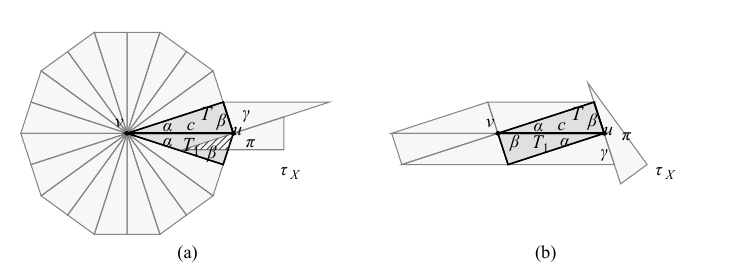}
\caption{Longest-side pairing}\label{fig:longest-pairing}
\end{minipage}\par\addvspace{\intextsep}

Consider the $22$ candidate vertex types listed in Appendix Table~\ref{tab:39}. Vertices containing $\beta$ and $\gamma$ have degree at least $4$ and $3$, respectively, and a degree-$4$ vertex containing $\beta$ must be of type $\tau_X$. If a pure $\alpha$ vertex of degree $20$ occurs, every incident tile $T$ has $q_\alpha=20$. Equation~(1) and $N=21$ give
\[
q_\beta+q_\gamma+\sum_{x\in\mathcal V_\pi(T)}(q(x)-2)=7.
\]
Since $q_\beta\ge4$, $q_\gamma\ge3$, and every subdivision vertex makes a positive contribution, we obtain $q_\beta=4$, $q_\gamma=3$, and no subdivision vertices on $T$. The $\beta$-corner of each incident tile is therefore a vertex of type $\tau_X$, so the preceding pairing count gives $20\le0$, a contradiction. All the remaining $21$ candidate vertex types satisfy
\[
4q(q-2)\le63r+12s+8t.
\]
Weighted summation using Lemma~\ref{lem:angle-incidence}\textup{(i)} and equation~\eqref{eq:angle-incidence} gives $N\le83/4<21$, excluding this prototile.

For the $15^\circ$-$45^\circ$-$120^\circ$ triangle, the bound $q(v)\le N=21$ from Lemma~\ref{lem:neighbors} excludes both vertices of degree $24$ and vertices of degree $22$.

The local argument for $N\le22$ in the proof of Proposition~\ref{prop:15} excludes $\tau_W$. All the remaining $23$ candidate vertex types satisfy
\[
20q(q-2)\le341r+135s-60t+240P.
\]
Weighted summation using Lemma~\ref{lem:angle-incidence}\textup{(i)}, equation~\eqref{eq:angle-incidence}, and inequality~\eqref{eq:parity-bound} gives $N\le104/5<21$. The two inequalities require $44$ exact substitutions in total; the slacks are given in Cases C23 and C24 of Appendix Table~\ref{app:rigid-slacks-combined}. Thus every prototile with $\alpha<30^\circ$ has been excluded.

For $\alpha\ge30^\circ$, every vertex satisfies $3\le q\le12$. The inequality $(q-3)(q-12)\le0$ gives
\[
\frac1q\le\frac{15-q}{36},
\]
with equality only when $q=3$ or $q=12$. For any tile $T$, write $m=m(T)$. Substituting the preceding inequality into the definition of curvature and using equation~(1) gives
\[
\kappa(T)\le1-\frac{m}{2}
+\frac{15m-(21+2m)}{36}
=\frac{5(3-m)}{36}\le0.
\tag{14}\label{eq:rigidity-curvature}
\]
By equation~\eqref{eq:rigidity-curvature}, $\kappa(T)=0$ requires $m=3$ and all three corner degrees to belong to $\{3,12\}$. Equation~(1) then gives a degree sum of $27$, so $\kappa(T)=0$ if and only if $T$ is of type $(3,12,12)$. We call such a tile a \emph{special tile}. If there were no special tile, the bounds $m\le21$ and $3\le q\le12$ would allow only finitely many negative curvature values. The curvature would therefore have a uniform negative upper bound, contradicting the average-curvature conclusion of Lemma~\ref{lem:curvature}. Hence at least one special tile exists.

At each of the two degree-$12$ vertices of a special tile, twelve angles of at least $30^\circ$ meet, so all these angles are $30^\circ$. The prototile is therefore a $30^\circ$-$30^\circ$-$120^\circ$ triangle. Denote its short and long side lengths by $a$ and $c$, respectively. Then
\[
c=\sqrt3\,a,\qquad a<c<2a.
\]
The vertex angle equation becomes
\[
r+s+4t+6\varepsilon=12,\qquad q=r+s+t+\varepsilon,
\]
whose nonnegative integer solutions restrict the vertex degrees to $\{3,4,6,7,9,12\}$. A degree-$3$ vertex consists of three $120^\circ$ angles.

A special tile has no subdivision vertices, and the prototile has no longer side. Hence the tile across its long side shares that entire side. Both endpoints of the shared side have degree $12$, contributing $20$ in total to equation~(1). Since $N=21$, the third corner of the neighboring tile has degree $3$, and the tile has no subdivision vertices. Hence it is also special.

Let $T_1$ be the tile across a short side of a special tile. Denote the degree-$3$ and degree-$12$ endpoints of the shared side by $p,u$, respectively, and the other $30^\circ$ corner of $T_1$ by $v$. At $p$, all three tiles contribute $120^\circ$ angles, and all sides issuing from $p$ are short sides of length $a$. They therefore match along entire sides, so both short sides of $T_1$ have no subdivision vertices. Along the long side $uv$, the tile on the opposite side contributes a $30^\circ$ angle at $u$, since $q(u)=12$, and its side issuing from $u$ along $uv$ has length $a$ or $c$. If its length is $c$, the long sides match along their entire lengths. Then $T_1$ has no subdivision vertices, and equation~(1) gives $q(v)=12$, so $T_1$ is special.

Suppose that the opposite side has length $a$, and let $z$ be its other endpoint. As Fig.~\ref{fig:short-edge-obstruction} shows, the remaining $60^\circ$ at $z$ must be filled by two $30^\circ$ angles, so $q(z)=4$. The next tile side from $z$ toward $v$ has length at least $a$; the comparison $|zv|=c-a<a$ in the figure therefore forces it to extend beyond $v$. The straight angles on both sides fill a full turn along $(z,v)$, leaving no further tiling vertex. Thus the boundary of $T_1$ contains exactly $u,p,z,v$, and equation~(1) gives $q(v)=10$, contrary to the allowed degree set. Hence the opposite side has length $c$, and the tile across a short side is also special.

\begin{center}
\begin{minipage}{\linewidth}
\centering
\includegraphics[width=70mm]{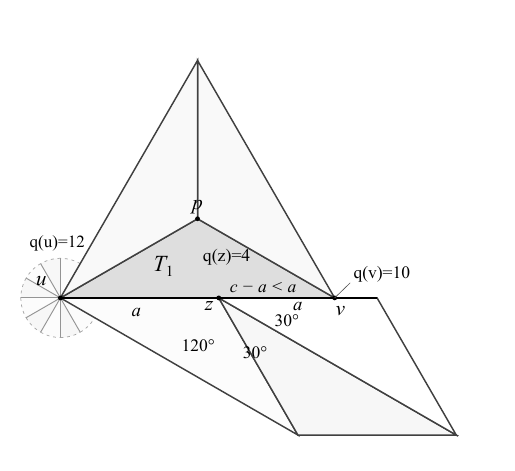}
\captionof{figure}{Local configuration under the assumed short-side match}
\label{fig:short-edge-obstruction}
\end{minipage}
\end{center}

Every tile sharing a side with a special tile is therefore special. By local finiteness, interior points of any two tiles can be joined by a polygonal path that avoids tiling vertices and crosses only finitely many tile sides transversely. The edge-adjacency graph, whose vertices are tiles and whose edges join pairs sharing a tiling edge, is therefore connected. Propagating along adjacency chains shows that every tile is special.

At every degree-$3$ vertex, three tiles match along their entire short sides and form an equilateral triangle of side length $c$. Each small tile has exactly one $120^\circ$ corner, so it belongs to a unique triple. These large triangles cover the plane with pairwise disjoint interiors and meet along entire, unsubdivided long sides, forming a side-to-side triangular tiling. The connectivity of the edge-adjacency graph of the small tiles implies that the large triangles are also connected through shared sides. Once one large triangle is fixed, the equilateral triangle across each side is uniquely determined by reflection in that side. Hence all the large triangles obtained along adjacency chains belong to the same standard triangular tiling. Since they cover the plane, they form the entire tiling. The common point of each triple of small tiles is equidistant from the three vertices of the large triangle, so it is its center. Thus the original tiling is similar to the central subdivision of the standard equilateral triangular tiling shown in Fig.~\ref{fig:construction}, proving Theorem~\ref{thm:main}.
\qedhere
\end{proof}

% Only the existing supplementary-information note and the requested AI declaration.
\textbf{Supplementary information.} The arithmetic verification code and input data are available at \url{https://github.com/shuise888/convex-monohedral-tilings}.

\section*{Declaration on the use of generative AI}
\PaperAIStatement\par

\phantomsection

\section*{Appendix}\label{start:appendix}\addcontentsline{toc}{section}{Appendix}

This appendix gives the exact slacks for the weighted inequalities and the successive vertex-type elimination records used in Sections~\ref{sec:triangles} and~\ref{sec:rigidity}.
\phantomsection\label{app:ordinary-checks}
For each rational-angle prototile, write $(\alpha,\beta,\gamma)=180^\circ(k_1,k_2,k_3)/\Lambda$ with positive integers $k_1,k_2,k_3$ and $\Lambda=k_1+k_2+k_3$. Choose a positive integer $L$ clearing the denominators in \eqref{eq:weighted-bound}. The slack is then $\Delta=L(ur+vs+wt+zP-q(q-2))$. By the concavity argument in Lemma~\ref{lem:weighted-screening}, only the two feasible endpoints of $s$ need to be checked in each group.
Tables~\ref{tab:8}--\ref{tab:38} group the vertex types by $(\varepsilon,t,P)$ and give the smallest and largest feasible values of $s$ and the corresponding integer slacks. When the endpoints coincide, only one slack is listed. Each table caption gives the total number of candidate vertex types at that stage; the Types column gives the number in the corresponding group. Deletion records use the order $(r,s,t,\varepsilon,q)$, and $\boldsymbol\mu=(\mu_1,\mu_2,\mu_3)$ records the minimum contribution $q-2$ for each of the three angles.

\begingroup
\tablefont
\setlength{\tabcolsep}{2pt}
\renewcommand{\arraystretch}{1}
\setlength{\LTpre}{2pt}
\setlength{\LTpost}{3pt}
\setlength{\LTcapwidth}{\linewidth}
\tablefont

\begin{longtable}{@{}>{\centering\arraybackslash}p{\dimexpr 0.135000\linewidth-2\tabcolsep\relax}>{\centering\arraybackslash}p{\dimexpr 0.105000\linewidth-2\tabcolsep\relax}>{\centering\arraybackslash}p{\dimexpr 0.190000\linewidth-2\tabcolsep\relax}>{\centering\arraybackslash}p{\dimexpr 0.070000\linewidth-2\tabcolsep\relax}>{\centering\arraybackslash}p{\dimexpr 0.135000\linewidth-2\tabcolsep\relax}>{\centering\arraybackslash}p{\dimexpr 0.105000\linewidth-2\tabcolsep\relax}>{\centering\arraybackslash}p{\dimexpr 0.190000\linewidth-2\tabcolsep\relax}>{\centering\arraybackslash}p{\dimexpr 0.070000\linewidth-2\tabcolsep\relax}@{}}
\caption[Endpoint checks, stage 1]{Endpoint checks, stage 1. Angles (in degrees) $(45/2, 135/2, 90)$; $K=\text{--}$; $(k_1,k_2,k_3)=(1, 3, 4)$, $\Lambda=8$; $L=1$, $L(u,v,w,z)=(14,3,2,1)$; $N\le 19$; 23 candidate vertex types.}\label{tab:8}\\[-2pt]
\toprule
$(\varepsilon,t,P)$ & $(s_-,s_+)$ & $\Delta_-,\Delta_+$ & Types & $(\varepsilon,t,P)$ & $(s_-,s_+)$ & $\Delta_-,\Delta_+$ & Types \\
\midrule
\endfirsthead
\caption[]{Continued}\\[-2pt]
\toprule
$(\varepsilon,t,P)$ & $(s_-,s_+)$ & $\Delta_-,\Delta_+$ & Types & $(\varepsilon,t,P)$ & $(s_-,s_+)$ & $\Delta_-,\Delta_+$ & Types \\
\midrule
\endhead
\midrule
\endfoot
\bottomrule
\endlastfoot
$(0, 0, 0)$ & $(0,4)$ & $0, 20$ & $3$ & $(0, 3, 0)$ & $(1,1)$ & $8$ & $1$ \\
$(0, 0, 1)$ & $(1,5)$ & $18, 6$ & $3$ & $(0, 3, 1)$ & $(0,0)$ & $28$ & $1$ \\
$(0, 1, 0)$ & $(1,3)$ & $32, 18$ & $2$ & $(0, 4, 0)$ & $(0,0)$ & $0$ & $1$ \\
$(0, 1, 1)$ & $(0,4)$ & $28, 0$ & $3$ & $(1, 0, -1)$ & $(0,2)$ & $48, 18$ & $3$ \\
$(0, 2, 0)$ & $(0,2)$ & $36, 14$ & $2$ & $(1, 1, -1)$ & $(0,1)$ & $33, 10$ & $2$ \\
$(0, 2, 1)$ & $(1,1)$ & $30$ & $1$ & $(1, 2, -1)$ & $(0,0)$ & $0$ & $1$ \\
\end{longtable}

\endgroup

\begingroup
\tablefont
\setlength{\tabcolsep}{2pt}
\renewcommand{\arraystretch}{1}
\setlength{\LTpre}{2pt}
\setlength{\LTpost}{3pt}
\setlength{\LTcapwidth}{\linewidth}
\tablefont

\begin{longtable}{@{}>{\centering\arraybackslash}p{\dimexpr 0.135000\linewidth-2\tabcolsep\relax}>{\centering\arraybackslash}p{\dimexpr 0.105000\linewidth-2\tabcolsep\relax}>{\centering\arraybackslash}p{\dimexpr 0.190000\linewidth-2\tabcolsep\relax}>{\centering\arraybackslash}p{\dimexpr 0.070000\linewidth-2\tabcolsep\relax}>{\centering\arraybackslash}p{\dimexpr 0.135000\linewidth-2\tabcolsep\relax}>{\centering\arraybackslash}p{\dimexpr 0.105000\linewidth-2\tabcolsep\relax}>{\centering\arraybackslash}p{\dimexpr 0.190000\linewidth-2\tabcolsep\relax}>{\centering\arraybackslash}p{\dimexpr 0.070000\linewidth-2\tabcolsep\relax}@{}}
\caption[Endpoint checks, stage 2]{Endpoint checks, stage 2. Angles (in degrees) $(90/7, 360/7, 810/7)$; $K=\text{--}$; $(k_1,k_2,k_3)=(1, 4, 9)$, $\Lambda=14$; $L=31$, $L(u,v,w,z)=(806,104,-305,357)$; $N\le 605/31$; 23 candidate vertex types.}\label{tab:9}\\[-2pt]
\toprule
$(\varepsilon,t,P)$ & $(s_-,s_+)$ & $\Delta_-,\Delta_+$ & Types & $(\varepsilon,t,P)$ & $(s_-,s_+)$ & $\Delta_-,\Delta_+$ & Types \\
\midrule
\endfirsthead
\caption[]{Continued}\\[-2pt]
\toprule
$(\varepsilon,t,P)$ & $(s_-,s_+)$ & $\Delta_-,\Delta_+$ & Types & $(\varepsilon,t,P)$ & $(s_-,s_+)$ & $\Delta_-,\Delta_+$ & Types \\
\midrule
\endhead
\midrule
\endfoot
\bottomrule
\endlastfoot
$(0, 0, 0)$ & $(0,6)$ & $0, 1368$ & $4$ & $(0, 2, 1)$ & $(1,1)$ & $2734$ & $1$ \\
$(0, 0, 1)$ & $(1,7)$ & $1980, 0$ & $4$ & $(0, 3, 1)$ & $(0,0)$ & $0$ & $1$ \\
$(0, 1, 0)$ & $(1,3)$ & $3984, 2580$ & $2$ & $(1, 0, -1)$ & $(0,3)$ & $4882, 823$ & $4$ \\
$(0, 1, 1)$ & $(0,4)$ & $4206, 1398$ & $3$ & $(1, 1, -1)$ & $(0,1)$ & $2283, 0$ & $2$ \\
$(0, 2, 0)$ & $(0,2)$ & $3730, 466$ & $2$ &  &  &  &  \\
\end{longtable}

\endgroup

\begingroup
\tablefont
\setlength{\tabcolsep}{2pt}
\renewcommand{\arraystretch}{1}
\setlength{\LTpre}{2pt}
\setlength{\LTpost}{3pt}
\setlength{\LTcapwidth}{\linewidth}
\tablefont

\begin{longtable}{@{}>{\centering\arraybackslash}p{\dimexpr 0.135000\linewidth-2\tabcolsep\relax}>{\centering\arraybackslash}p{\dimexpr 0.105000\linewidth-2\tabcolsep\relax}>{\centering\arraybackslash}p{\dimexpr 0.190000\linewidth-2\tabcolsep\relax}>{\centering\arraybackslash}p{\dimexpr 0.070000\linewidth-2\tabcolsep\relax}>{\centering\arraybackslash}p{\dimexpr 0.135000\linewidth-2\tabcolsep\relax}>{\centering\arraybackslash}p{\dimexpr 0.105000\linewidth-2\tabcolsep\relax}>{\centering\arraybackslash}p{\dimexpr 0.190000\linewidth-2\tabcolsep\relax}>{\centering\arraybackslash}p{\dimexpr 0.070000\linewidth-2\tabcolsep\relax}@{}}
\caption[Endpoint checks, stage 3]{Endpoint checks, stage 3. Angles (in degrees) $(180/11, 540/11, 1260/11)$; $K=\text{--}$; $(k_1,k_2,k_3)=(1, 3, 7)$, $\Lambda=11$; $L=11$, $L(u,v,w,z)=(220,48,-78,102)$; $N\le 190/11$; 24 candidate vertex types.}\label{tab:10}\\[-2pt]
\toprule
$(\varepsilon,t,P)$ & $(s_-,s_+)$ & $\Delta_-,\Delta_+$ & Types & $(\varepsilon,t,P)$ & $(s_-,s_+)$ & $\Delta_-,\Delta_+$ & Types \\
\midrule
\endfirsthead
\caption[]{Continued}\\[-2pt]
\toprule
$(\varepsilon,t,P)$ & $(s_-,s_+)$ & $\Delta_-,\Delta_+$ & Types & $(\varepsilon,t,P)$ & $(s_-,s_+)$ & $\Delta_-,\Delta_+$ & Types \\
\midrule
\endhead
\midrule
\endfoot
\bottomrule
\endlastfoot
$(0, 0, 0)$ & $(0,6)$ & $0, 288$ & $4$ & $(0, 2, 1)$ & $(1,1)$ & $566$ & $1$ \\
$(0, 0, 1)$ & $(1,7)$ & $370, 130$ & $4$ & $(0, 3, 1)$ & $(0,0)$ & $0$ & $1$ \\
$(0, 1, 1)$ & $(0,5)$ & $860, 0$ & $6$ & $(1, 0, -1)$ & $(0,3)$ & $998, 218$ & $4$ \\
$(0, 2, 0)$ & $(0,2)$ & $724, 116$ & $2$ & $(1, 1, -1)$ & $(0,1)$ & $436, 0$ & $2$ \\
\end{longtable}

\endgroup

\begingroup
\tablefont
\setlength{\tabcolsep}{2pt}
\renewcommand{\arraystretch}{1}
\setlength{\LTpre}{2pt}
\setlength{\LTpost}{3pt}
\setlength{\LTcapwidth}{\linewidth}
\tablefont

\begin{longtable}{@{}>{\centering\arraybackslash}p{\dimexpr 0.135000\linewidth-2\tabcolsep\relax}>{\centering\arraybackslash}p{\dimexpr 0.105000\linewidth-2\tabcolsep\relax}>{\centering\arraybackslash}p{\dimexpr 0.190000\linewidth-2\tabcolsep\relax}>{\centering\arraybackslash}p{\dimexpr 0.070000\linewidth-2\tabcolsep\relax}>{\centering\arraybackslash}p{\dimexpr 0.135000\linewidth-2\tabcolsep\relax}>{\centering\arraybackslash}p{\dimexpr 0.105000\linewidth-2\tabcolsep\relax}>{\centering\arraybackslash}p{\dimexpr 0.190000\linewidth-2\tabcolsep\relax}>{\centering\arraybackslash}p{\dimexpr 0.070000\linewidth-2\tabcolsep\relax}@{}}
\caption[Endpoint checks, stage 4]{Endpoint checks, stage 4. Angles (in degrees) $(45/2, 45, 225/2)$; $K=\text{--}$; $(k_1,k_2,k_3)=(1, 2, 5)$, $\Lambda=8$; $L=2$, $L(u,v,w,z)=(28,12,-9,15)$; $N\le 31/2$; 27 candidate vertex types.}\label{tab:11}\\[-2pt]
\toprule
$(\varepsilon,t,P)$ & $(s_-,s_+)$ & $\Delta_-,\Delta_+$ & Types & $(\varepsilon,t,P)$ & $(s_-,s_+)$ & $\Delta_-,\Delta_+$ & Types \\
\midrule
\endfirsthead
\caption[]{Continued}\\[-2pt]
\toprule
$(\varepsilon,t,P)$ & $(s_-,s_+)$ & $\Delta_-,\Delta_+$ & Types & $(\varepsilon,t,P)$ & $(s_-,s_+)$ & $\Delta_-,\Delta_+$ & Types \\
\midrule
\endhead
\midrule
\endfoot
\bottomrule
\endlastfoot
$(0, 0, 0)$ & $(0,8)$ & $0, 0$ & $5$ & $(0, 2, 1)$ & $(1,3)$ & $51, 3$ & $2$ \\
$(0, 0, 1)$ & $(1,7)$ & $29, 29$ & $4$ & $(0, 3, 1)$ & $(0,0)$ & $0$ & $1$ \\
$(0, 1, 0)$ & $(1,5)$ & $57, 9$ & $3$ & $(1, 0, -1)$ & $(0,4)$ & $83, 3$ & $5$ \\
$(0, 1, 1)$ & $(0,4)$ & $74, 42$ & $3$ & $(1, 1, -1)$ & $(0,1)$ & $30, 0$ & $2$ \\
$(0, 2, 0)$ & $(0,2)$ & $54, 14$ & $2$ &  &  &  &  \\
\end{longtable}

\endgroup

\begingroup
\tablefont
\setlength{\tabcolsep}{2pt}
\renewcommand{\arraystretch}{1}
\setlength{\LTpre}{2pt}
\setlength{\LTpost}{3pt}
\setlength{\LTcapwidth}{\linewidth}
\tablefont

\begin{longtable}{@{}>{\centering\arraybackslash}p{\dimexpr 0.135000\linewidth-2\tabcolsep\relax}>{\centering\arraybackslash}p{\dimexpr 0.105000\linewidth-2\tabcolsep\relax}>{\centering\arraybackslash}p{\dimexpr 0.190000\linewidth-2\tabcolsep\relax}>{\centering\arraybackslash}p{\dimexpr 0.070000\linewidth-2\tabcolsep\relax}>{\centering\arraybackslash}p{\dimexpr 0.135000\linewidth-2\tabcolsep\relax}>{\centering\arraybackslash}p{\dimexpr 0.105000\linewidth-2\tabcolsep\relax}>{\centering\arraybackslash}p{\dimexpr 0.190000\linewidth-2\tabcolsep\relax}>{\centering\arraybackslash}p{\dimexpr 0.070000\linewidth-2\tabcolsep\relax}@{}}
\caption[Endpoint checks, stage 5]{Endpoint checks, stage 5. Angles (in degrees) $(360/13, 540/13, 1440/13)$; $K=\text{--}$; $(k_1,k_2,k_3)=(2, 3, 8)$, $\Lambda=13$; $L=13$, $L(u,v,w,z)=(140,80,-38,78)$; $N\le 14$; 15 candidate vertex types.}\label{tab:12}\\[-2pt]
\toprule
$(\varepsilon,t,P)$ & $(s_-,s_+)$ & $\Delta_-,\Delta_+$ & Types & $(\varepsilon,t,P)$ & $(s_-,s_+)$ & $\Delta_-,\Delta_+$ & Types \\
\midrule
\endfirsthead
\caption[]{Continued}\\[-2pt]
\toprule
$(\varepsilon,t,P)$ & $(s_-,s_+)$ & $\Delta_-,\Delta_+$ & Types & $(\varepsilon,t,P)$ & $(s_-,s_+)$ & $\Delta_-,\Delta_+$ & Types \\
\midrule
\endhead
\midrule
\endfoot
\bottomrule
\endlastfoot
$(0, 0, 0)$ & $(2,6)$ & $0, 0$ & $2$ & $(0, 2, 1)$ & $(0,0)$ & $247$ & $1$ \\
$(0, 0, 1)$ & $(0,8)$ & $39, 39$ & $3$ & $(0, 3, 1)$ & $(0,0)$ & $0$ & $1$ \\
$(0, 1, 1)$ & $(0,6)$ & $260, 65$ & $4$ & $(1, 0, -1)$ & $(1,3)$ & $247, 130$ & $2$ \\
$(0, 2, 0)$ & $(2,2)$ & $52$ & $1$ & $(1, 1, -1)$ & $(1,1)$ & $0$ & $1$ \\
\end{longtable}

\endgroup

\begingroup
\tablefont
\setlength{\tabcolsep}{2pt}
\renewcommand{\arraystretch}{1}
\setlength{\LTpre}{2pt}
\setlength{\LTpost}{3pt}
\setlength{\LTcapwidth}{\linewidth}
\tablefont

\begin{longtable}{@{}>{\centering\arraybackslash}p{\dimexpr 0.135000\linewidth-2\tabcolsep\relax}>{\centering\arraybackslash}p{\dimexpr 0.105000\linewidth-2\tabcolsep\relax}>{\centering\arraybackslash}p{\dimexpr 0.190000\linewidth-2\tabcolsep\relax}>{\centering\arraybackslash}p{\dimexpr 0.070000\linewidth-2\tabcolsep\relax}>{\centering\arraybackslash}p{\dimexpr 0.135000\linewidth-2\tabcolsep\relax}>{\centering\arraybackslash}p{\dimexpr 0.105000\linewidth-2\tabcolsep\relax}>{\centering\arraybackslash}p{\dimexpr 0.190000\linewidth-2\tabcolsep\relax}>{\centering\arraybackslash}p{\dimexpr 0.070000\linewidth-2\tabcolsep\relax}@{}}
\caption[Endpoint checks, stage 6]{Endpoint checks, stage 6. Angles (in degrees) $(12, 72, 96)$; $K=\text{--}$; $(k_1,k_2,k_3)=(1, 6, 8)$, $\Lambda=15$; $L=11$, $L(u,v,w,z)=(308,-4,-31,185)$; $N\le 273/11$; 20 candidate vertex types.}\label{tab:13}\\[-2pt]
\toprule
$(\varepsilon,t,P)$ & $(s_-,s_+)$ & $\Delta_-,\Delta_+$ & Types & $(\varepsilon,t,P)$ & $(s_-,s_+)$ & $\Delta_-,\Delta_+$ & Types \\
\midrule
\endfirsthead
\caption[]{Continued}\\[-2pt]
\toprule
$(\varepsilon,t,P)$ & $(s_-,s_+)$ & $\Delta_-,\Delta_+$ & Types & $(\varepsilon,t,P)$ & $(s_-,s_+)$ & $\Delta_-,\Delta_+$ & Types \\
\midrule
\endhead
\midrule
\endfoot
\bottomrule
\endlastfoot
$(0, 0, 0)$ & $(0,4)$ & $0, 952$ & $3$ & $(0, 2, 1)$ & $(1,1)$ & $1494$ & $1$ \\
$(0, 0, 1)$ & $(1,5)$ & $1248, 0$ & $3$ & $(0, 3, 1)$ & $(0,1)$ & $1247, 0$ & $2$ \\
$(0, 1, 1)$ & $(0,3)$ & $1617, 846$ & $4$ & $(1, 0, -1)$ & $(0,2)$ & $1971, 467$ & $3$ \\
$(0, 2, 0)$ & $(0,2)$ & $1786, 282$ & $2$ & $(1, 1, -1)$ & $(0,1)$ & $1247, 0$ & $2$ \\
\end{longtable}

\endgroup

\begingroup
\tablefont
\setlength{\tabcolsep}{2pt}
\renewcommand{\arraystretch}{1}
\setlength{\LTpre}{2pt}
\setlength{\LTpost}{3pt}
\setlength{\LTcapwidth}{\linewidth}
\tablefont

\begin{longtable}{@{}>{\centering\arraybackslash}p{\dimexpr 0.135000\linewidth-2\tabcolsep\relax}>{\centering\arraybackslash}p{\dimexpr 0.105000\linewidth-2\tabcolsep\relax}>{\centering\arraybackslash}p{\dimexpr 0.190000\linewidth-2\tabcolsep\relax}>{\centering\arraybackslash}p{\dimexpr 0.070000\linewidth-2\tabcolsep\relax}>{\centering\arraybackslash}p{\dimexpr 0.135000\linewidth-2\tabcolsep\relax}>{\centering\arraybackslash}p{\dimexpr 0.105000\linewidth-2\tabcolsep\relax}>{\centering\arraybackslash}p{\dimexpr 0.190000\linewidth-2\tabcolsep\relax}>{\centering\arraybackslash}p{\dimexpr 0.070000\linewidth-2\tabcolsep\relax}@{}}
\caption[Endpoint checks, stage 7]{Endpoint checks, stage 7. Angles (in degrees) $(12, 72, 96)$; $K=24$; $(k_1,k_2,k_3)=(1, 6, 8)$, $\Lambda=15$; $L=195$, $L(u,v,w,z)=(3884,144,-263,2205)$; $N\le 251/13$; 17 candidate vertex types. $\boldsymbol\mu=(2, 2, 2)$. The neighbor-number constraint excludes the following 3 vertex types: [(22, 0, 1, 0, 23), (24, 1, 0, 0, 25), (30, 0, 0, 0, 30)]. $\boldsymbol\mu=(2, 2, 2)$. Additional vertex types deleted in this round: 0; elimination terminates.}\label{tab:14}\\[-2pt]
\toprule
$(\varepsilon,t,P)$ & $(s_-,s_+)$ & $\Delta_-,\Delta_+$ & Types & $(\varepsilon,t,P)$ & $(s_-,s_+)$ & $\Delta_-,\Delta_+$ & Types \\
\midrule
\endfirsthead
\caption[]{Continued}\\[-2pt]
\toprule
$(\varepsilon,t,P)$ & $(s_-,s_+)$ & $\Delta_-,\Delta_+$ & Types & $(\varepsilon,t,P)$ & $(s_-,s_+)$ & $\Delta_-,\Delta_+$ & Types \\
\midrule
\endhead
\midrule
\endfoot
\bottomrule
\endlastfoot
$(0, 0, 0)$ & $(2,4)$ & $0, 8280$ & $2$ & $(0, 2, 1)$ & $(1,1)$ & $13590$ & $1$ \\
$(0, 0, 1)$ & $(3,5)$ & $11220, 0$ & $2$ & $(0, 3, 1)$ & $(0,1)$ & $12435, 0$ & $2$ \\
$(0, 1, 1)$ & $(1,3)$ & $8070, 8550$ & $3$ & $(1, 0, -1)$ & $(0,2)$ & $12375, 5055$ & $3$ \\
$(0, 2, 0)$ & $(0,2)$ & $10170, 2850$ & $2$ & $(1, 1, -1)$ & $(0,1)$ & $12435, 0$ & $2$ \\
\end{longtable}

\endgroup

\begingroup
\tablefont
\setlength{\tabcolsep}{2pt}
\renewcommand{\arraystretch}{1}
\setlength{\LTpre}{2pt}
\setlength{\LTpost}{3pt}
\setlength{\LTcapwidth}{\linewidth}
\tablefont

\begin{longtable}{@{}>{\centering\arraybackslash}p{\dimexpr 0.135000\linewidth-2\tabcolsep\relax}>{\centering\arraybackslash}p{\dimexpr 0.105000\linewidth-2\tabcolsep\relax}>{\centering\arraybackslash}p{\dimexpr 0.190000\linewidth-2\tabcolsep\relax}>{\centering\arraybackslash}p{\dimexpr 0.070000\linewidth-2\tabcolsep\relax}>{\centering\arraybackslash}p{\dimexpr 0.135000\linewidth-2\tabcolsep\relax}>{\centering\arraybackslash}p{\dimexpr 0.105000\linewidth-2\tabcolsep\relax}>{\centering\arraybackslash}p{\dimexpr 0.190000\linewidth-2\tabcolsep\relax}>{\centering\arraybackslash}p{\dimexpr 0.070000\linewidth-2\tabcolsep\relax}@{}}
\caption[Endpoint checks, stage 8]{Endpoint checks, stage 8. Angles (in degrees) $(180/13, 900/13, 1260/13)$; $K=\text{--}$; $(k_1,k_2,k_3)=(1, 5, 7)$, $\Lambda=13$; $L=11$, $L(u,v,w,z)=(264,-28,-8,140)$; $N\le 228/11$; 20 candidate vertex types.}\label{tab:15}\\[-2pt]
\toprule
$(\varepsilon,t,P)$ & $(s_-,s_+)$ & $\Delta_-,\Delta_+$ & Types & $(\varepsilon,t,P)$ & $(s_-,s_+)$ & $\Delta_-,\Delta_+$ & Types \\
\midrule
\endfirsthead
\caption[]{Continued}\\[-2pt]
\toprule
$(\varepsilon,t,P)$ & $(s_-,s_+)$ & $\Delta_-,\Delta_+$ & Types & $(\varepsilon,t,P)$ & $(s_-,s_+)$ & $\Delta_-,\Delta_+$ & Types \\
\midrule
\endhead
\midrule
\endfoot
\bottomrule
\endlastfoot
$(0, 0, 0)$ & $(0,4)$ & $0, 592$ & $3$ & $(0, 2, 1)$ & $(1,1)$ & $1064$ & $1$ \\
$(0, 0, 1)$ & $(1,5)$ & $816, 0$ & $3$ & $(0, 3, 1)$ & $(0,1)$ & $908, 0$ & $2$ \\
$(0, 1, 1)$ & $(0,3)$ & $1188, 576$ & $4$ & $(1, 0, -1)$ & $(0,2)$ & $1444, 332$ & $3$ \\
$(0, 2, 0)$ & $(0,2)$ & $1304, 192$ & $2$ & $(1, 1, -1)$ & $(0,1)$ & $908, 0$ & $2$ \\
\end{longtable}

\endgroup

\begingroup
\tablefont
\setlength{\tabcolsep}{2pt}
\renewcommand{\arraystretch}{1}
\setlength{\LTpre}{2pt}
\setlength{\LTpost}{3pt}
\setlength{\LTcapwidth}{\linewidth}
\tablefont

\begin{longtable}{@{}>{\centering\arraybackslash}p{\dimexpr 0.135000\linewidth-2\tabcolsep\relax}>{\centering\arraybackslash}p{\dimexpr 0.105000\linewidth-2\tabcolsep\relax}>{\centering\arraybackslash}p{\dimexpr 0.190000\linewidth-2\tabcolsep\relax}>{\centering\arraybackslash}p{\dimexpr 0.070000\linewidth-2\tabcolsep\relax}>{\centering\arraybackslash}p{\dimexpr 0.135000\linewidth-2\tabcolsep\relax}>{\centering\arraybackslash}p{\dimexpr 0.105000\linewidth-2\tabcolsep\relax}>{\centering\arraybackslash}p{\dimexpr 0.190000\linewidth-2\tabcolsep\relax}>{\centering\arraybackslash}p{\dimexpr 0.070000\linewidth-2\tabcolsep\relax}@{}}
\caption[Endpoint checks, stage 9]{Endpoint checks, stage 9. Angles (in degrees) $(180/11, 720/11, 1080/11)$; $K=\text{--}$; $(k_1,k_2,k_3)=(1, 4, 6)$, $\Lambda=11$; $L=8$, $L(u,v,w,z)=(160,10,-13,93)$; $N\le 157/8$; 21 candidate vertex types.}\label{tab:16}\\[-2pt]
\toprule
$(\varepsilon,t,P)$ & $(s_-,s_+)$ & $\Delta_-,\Delta_+$ & Types & $(\varepsilon,t,P)$ & $(s_-,s_+)$ & $\Delta_-,\Delta_+$ & Types \\
\midrule
\endfirsthead
\caption[]{Continued}\\[-2pt]
\toprule
$(\varepsilon,t,P)$ & $(s_-,s_+)$ & $\Delta_-,\Delta_+$ & Types & $(\varepsilon,t,P)$ & $(s_-,s_+)$ & $\Delta_-,\Delta_+$ & Types \\
\midrule
\endhead
\midrule
\endfoot
\bottomrule
\endlastfoot
$(0, 0, 0)$ & $(0,4)$ & $0, 360$ & $3$ & $(0, 2, 1)$ & $(1,1)$ & $533$ & $1$ \\
$(0, 0, 1)$ & $(1,5)$ & $399, 183$ & $3$ & $(0, 3, 1)$ & $(0,1)$ & $414, 0$ & $2$ \\
$(0, 1, 1)$ & $(0,4)$ & $600, 0$ & $5$ & $(1, 0, -1)$ & $(0,2)$ & $707, 215$ & $3$ \\
$(0, 2, 0)$ & $(0,2)$ & $614, 122$ & $2$ & $(1, 1, -1)$ & $(0,1)$ & $414, 0$ & $2$ \\
\end{longtable}

\endgroup

\begingroup
\tablefont
\setlength{\tabcolsep}{2pt}
\renewcommand{\arraystretch}{1}
\setlength{\LTpre}{2pt}
\setlength{\LTpost}{3pt}
\setlength{\LTcapwidth}{\linewidth}
\tablefont

\begin{longtable}{@{}>{\centering\arraybackslash}p{\dimexpr 0.135000\linewidth-2\tabcolsep\relax}>{\centering\arraybackslash}p{\dimexpr 0.105000\linewidth-2\tabcolsep\relax}>{\centering\arraybackslash}p{\dimexpr 0.190000\linewidth-2\tabcolsep\relax}>{\centering\arraybackslash}p{\dimexpr 0.070000\linewidth-2\tabcolsep\relax}>{\centering\arraybackslash}p{\dimexpr 0.135000\linewidth-2\tabcolsep\relax}>{\centering\arraybackslash}p{\dimexpr 0.105000\linewidth-2\tabcolsep\relax}>{\centering\arraybackslash}p{\dimexpr 0.190000\linewidth-2\tabcolsep\relax}>{\centering\arraybackslash}p{\dimexpr 0.070000\linewidth-2\tabcolsep\relax}@{}}
\caption[Endpoint checks, stage 10]{Endpoint checks, stage 10. Angles (in degrees) $(20, 60, 100)$; $K=\text{--}$; $(k_1,k_2,k_3)=(1, 3, 5)$, $\Lambda=9$; $L=5$, $L(u,v,w,z)=(80,32,-16,56)$; $N\le 96/5$; 23 candidate vertex types.}\label{tab:17}\\[-2pt]
\toprule
$(\varepsilon,t,P)$ & $(s_-,s_+)$ & $\Delta_-,\Delta_+$ & Types & $(\varepsilon,t,P)$ & $(s_-,s_+)$ & $\Delta_-,\Delta_+$ & Types \\
\midrule
\endfirsthead
\caption[]{Continued}\\[-2pt]
\toprule
$(\varepsilon,t,P)$ & $(s_-,s_+)$ & $\Delta_-,\Delta_+$ & Types & $(\varepsilon,t,P)$ & $(s_-,s_+)$ & $\Delta_-,\Delta_+$ & Types \\
\midrule
\endhead
\midrule
\endfoot
\bottomrule
\endlastfoot
$(0, 0, 0)$ & $(0,6)$ & $0, 72$ & $4$ & $(0, 2, 1)$ & $(1,1)$ & $216$ & $1$ \\
$(0, 0, 1)$ & $(1,5)$ & $168, 216$ & $3$ & $(0, 3, 1)$ & $(0,1)$ & $128, 0$ & $2$ \\
$(0, 1, 1)$ & $(0,4)$ & $240, 128$ & $5$ & $(1, 0, -1)$ & $(0,3)$ & $264, 0$ & $4$ \\
$(0, 2, 0)$ & $(0,2)$ & $208, 72$ & $2$ & $(1, 1, -1)$ & $(0,1)$ & $128, 0$ & $2$ \\
\end{longtable}

\endgroup

\begingroup
\tablefont
\setlength{\tabcolsep}{2pt}
\renewcommand{\arraystretch}{1}
\setlength{\LTpre}{2pt}
\setlength{\LTpost}{3pt}
\setlength{\LTcapwidth}{\linewidth}
\tablefont

\begin{longtable}{@{}>{\centering\arraybackslash}p{\dimexpr 0.135000\linewidth-2\tabcolsep\relax}>{\centering\arraybackslash}p{\dimexpr 0.105000\linewidth-2\tabcolsep\relax}>{\centering\arraybackslash}p{\dimexpr 0.190000\linewidth-2\tabcolsep\relax}>{\centering\arraybackslash}p{\dimexpr 0.070000\linewidth-2\tabcolsep\relax}>{\centering\arraybackslash}p{\dimexpr 0.135000\linewidth-2\tabcolsep\relax}>{\centering\arraybackslash}p{\dimexpr 0.105000\linewidth-2\tabcolsep\relax}>{\centering\arraybackslash}p{\dimexpr 0.190000\linewidth-2\tabcolsep\relax}>{\centering\arraybackslash}p{\dimexpr 0.070000\linewidth-2\tabcolsep\relax}@{}}
\caption[Endpoint checks, stage 11]{Endpoint checks, stage 11. Angles (in degrees) $(45/2, 225/4, 405/4)$; $K=\text{--}$; $(k_1,k_2,k_3)=(2, 5, 9)$, $\Lambda=16$; $L=13$, $L(u,v,w,z)=(182,29,-8,99)$; $N\le 203/13$; 12 candidate vertex types.}\label{tab:18}\\[-2pt]
\toprule
$(\varepsilon,t,P)$ & $(s_-,s_+)$ & $\Delta_-,\Delta_+$ & Types & $(\varepsilon,t,P)$ & $(s_-,s_+)$ & $\Delta_-,\Delta_+$ & Types \\
\midrule
\endfirsthead
\caption[]{Continued}\\[-2pt]
\toprule
$(\varepsilon,t,P)$ & $(s_-,s_+)$ & $\Delta_-,\Delta_+$ & Types & $(\varepsilon,t,P)$ & $(s_-,s_+)$ & $\Delta_-,\Delta_+$ & Types \\
\midrule
\endhead
\midrule
\endfoot
\bottomrule
\endlastfoot
$(0, 0, 0)$ & $(0,4)$ & $0, 168$ & $2$ & $(0, 2, 1)$ & $(0,0)$ & $538$ & $1$ \\
$(0, 0, 1)$ & $(2,6)$ & $300, 0$ & $2$ & $(0, 3, 1)$ & $(1,1)$ & $0$ & $1$ \\
$(0, 1, 0)$ & $(1,1)$ & $372$ & $1$ & $(1, 0, -1)$ & $(0,2)$ & $538, 193$ & $2$ \\
$(0, 1, 1)$ & $(3,3)$ & $282$ & $1$ & $(1, 1, -1)$ & $(1,1)$ & $0$ & $1$ \\
$(0, 2, 0)$ & $(2,2)$ & $94$ & $1$ &  &  &  &  \\
\end{longtable}

\endgroup

\begingroup
\tablefont
\setlength{\tabcolsep}{2pt}
\renewcommand{\arraystretch}{1}
\setlength{\LTpre}{2pt}
\setlength{\LTpost}{3pt}
\setlength{\LTcapwidth}{\linewidth}
\tablefont

\begin{longtable}{@{}>{\centering\arraybackslash}p{\dimexpr 0.135000\linewidth-2\tabcolsep\relax}>{\centering\arraybackslash}p{\dimexpr 0.105000\linewidth-2\tabcolsep\relax}>{\centering\arraybackslash}p{\dimexpr 0.190000\linewidth-2\tabcolsep\relax}>{\centering\arraybackslash}p{\dimexpr 0.070000\linewidth-2\tabcolsep\relax}>{\centering\arraybackslash}p{\dimexpr 0.135000\linewidth-2\tabcolsep\relax}>{\centering\arraybackslash}p{\dimexpr 0.105000\linewidth-2\tabcolsep\relax}>{\centering\arraybackslash}p{\dimexpr 0.190000\linewidth-2\tabcolsep\relax}>{\centering\arraybackslash}p{\dimexpr 0.070000\linewidth-2\tabcolsep\relax}@{}}
\caption[Endpoint checks, stage 12]{Endpoint checks, stage 12. Angles (in degrees) $(180/7, 360/7, 720/7)$; $K=\text{--}$; $(k_1,k_2,k_3)=(1, 2, 4)$, $\Lambda=7$; $L=1$, $L(u,v,w,z)=(12,4,-1,7)$; $N\le 15$; 26 candidate vertex types.}\label{tab:19}\\[-2pt]
\toprule
$(\varepsilon,t,P)$ & $(s_-,s_+)$ & $\Delta_-,\Delta_+$ & Types & $(\varepsilon,t,P)$ & $(s_-,s_+)$ & $\Delta_-,\Delta_+$ & Types \\
\midrule
\endfirsthead
\caption[]{Continued}\\[-2pt]
\toprule
$(\varepsilon,t,P)$ & $(s_-,s_+)$ & $\Delta_-,\Delta_+$ & Types & $(\varepsilon,t,P)$ & $(s_-,s_+)$ & $\Delta_-,\Delta_+$ & Types \\
\midrule
\endhead
\midrule
\endfoot
\bottomrule
\endlastfoot
$(0, 0, 0)$ & $(0,6)$ & $0, 0$ & $4$ & $(0, 2, 1)$ & $(1,3)$ & $22, 2$ & $2$ \\
$(0, 0, 1)$ & $(1,7)$ & $12, 0$ & $4$ & $(0, 3, 1)$ & $(0,1)$ & $13, 0$ & $2$ \\
$(0, 1, 1)$ & $(0,5)$ & $27, 2$ & $6$ & $(1, 0, -1)$ & $(0,3)$ & $29, 2$ & $4$ \\
$(0, 2, 0)$ & $(0,2)$ & $22, 6$ & $2$ & $(1, 1, -1)$ & $(0,1)$ & $13, 0$ & $2$ \\
\end{longtable}

\endgroup

\begingroup
\tablefont
\setlength{\tabcolsep}{2pt}
\renewcommand{\arraystretch}{1}
\setlength{\LTpre}{2pt}
\setlength{\LTpost}{3pt}
\setlength{\LTcapwidth}{\linewidth}
\tablefont

\begin{longtable}{@{}>{\centering\arraybackslash}p{\dimexpr 0.135000\linewidth-2\tabcolsep\relax}>{\centering\arraybackslash}p{\dimexpr 0.105000\linewidth-2\tabcolsep\relax}>{\centering\arraybackslash}p{\dimexpr 0.190000\linewidth-2\tabcolsep\relax}>{\centering\arraybackslash}p{\dimexpr 0.070000\linewidth-2\tabcolsep\relax}>{\centering\arraybackslash}p{\dimexpr 0.135000\linewidth-2\tabcolsep\relax}>{\centering\arraybackslash}p{\dimexpr 0.105000\linewidth-2\tabcolsep\relax}>{\centering\arraybackslash}p{\dimexpr 0.190000\linewidth-2\tabcolsep\relax}>{\centering\arraybackslash}p{\dimexpr 0.070000\linewidth-2\tabcolsep\relax}@{}}
\caption[Endpoint checks, stage 13]{Endpoint checks, stage 13. Angles (in degrees) $(540/19, 900/19, 1980/19)$; $K=\text{--}$; $(k_1,k_2,k_3)=(3, 5, 11)$, $\Lambda=19$; $L=19$, $L(u,v,w,z)=(188,98,-20,114)$; $N\le 14$; 9 candidate vertex types.}\label{tab:20}\\[-2pt]
\toprule
$(\varepsilon,t,P)$ & $(s_-,s_+)$ & $\Delta_-,\Delta_+$ & Types & $(\varepsilon,t,P)$ & $(s_-,s_+)$ & $\Delta_-,\Delta_+$ & Types \\
\midrule
\endfirsthead
\caption[]{Continued}\\[-2pt]
\toprule
$(\varepsilon,t,P)$ & $(s_-,s_+)$ & $\Delta_-,\Delta_+$ & Types & $(\varepsilon,t,P)$ & $(s_-,s_+)$ & $\Delta_-,\Delta_+$ & Types \\
\midrule
\endhead
\midrule
\endfoot
\bottomrule
\endlastfoot
$(0, 0, 0)$ & $(4,4)$ & $0$ & $1$ & $(0, 3, 1)$ & $(1,1)$ & $0$ & $1$ \\
$(0, 0, 1)$ & $(1,7)$ & $0, 76$ & $2$ & $(1, 0, -1)$ & $(2,2)$ & $190$ & $1$ \\
$(0, 1, 1)$ & $(0,3)$ & $266, 228$ & $2$ & $(1, 1, -1)$ & $(1,1)$ & $0$ & $1$ \\
$(0, 2, 0)$ & $(2,2)$ & $76$ & $1$ &  &  &  &  \\
\end{longtable}

\endgroup

\begingroup
\tablefont
\setlength{\tabcolsep}{2pt}
\renewcommand{\arraystretch}{1}
\setlength{\LTpre}{2pt}
\setlength{\LTpost}{3pt}
\setlength{\LTcapwidth}{\linewidth}
\tablefont

\begin{longtable}{@{}>{\centering\arraybackslash}p{\dimexpr 0.135000\linewidth-2\tabcolsep\relax}>{\centering\arraybackslash}p{\dimexpr 0.105000\linewidth-2\tabcolsep\relax}>{\centering\arraybackslash}p{\dimexpr 0.190000\linewidth-2\tabcolsep\relax}>{\centering\arraybackslash}p{\dimexpr 0.070000\linewidth-2\tabcolsep\relax}>{\centering\arraybackslash}p{\dimexpr 0.135000\linewidth-2\tabcolsep\relax}>{\centering\arraybackslash}p{\dimexpr 0.105000\linewidth-2\tabcolsep\relax}>{\centering\arraybackslash}p{\dimexpr 0.190000\linewidth-2\tabcolsep\relax}>{\centering\arraybackslash}p{\dimexpr 0.070000\linewidth-2\tabcolsep\relax}@{}}
\caption[Endpoint checks, stage 14]{Endpoint checks, stage 14. Angles (in degrees) $(60/7, 360/7, 120)$; $K=\text{--}$; $(k_1,k_2,k_3)=(1, 6, 14)$, $\Lambda=21$; $L=31$, $L(u,v,w,z)=(1240,41,-235,798)$; $N\le 1046/31$; 23 candidate vertex types.}\label{tab:21}\\[-2pt]
\toprule
$(\varepsilon,t,P)$ & $(s_-,s_+)$ & $\Delta_-,\Delta_+$ & Types & $(\varepsilon,t,P)$ & $(s_-,s_+)$ & $\Delta_-,\Delta_+$ & Types \\
\midrule
\endfirsthead
\caption[]{Continued}\\[-2pt]
\toprule
$(\varepsilon,t,P)$ & $(s_-,s_+)$ & $\Delta_-,\Delta_+$ & Types & $(\varepsilon,t,P)$ & $(s_-,s_+)$ & $\Delta_-,\Delta_+$ & Types \\
\midrule
\endhead
\midrule
\endfoot
\bottomrule
\endlastfoot
$(0, 0, 0)$ & $(0,6)$ & $0, 3966$ & $4$ & $(0, 2, 1)$ & $(1,1)$ & $7220$ & $1$ \\
$(0, 0, 1)$ & $(1,7)$ & $5334, 0$ & $4$ & $(0, 3, 1)$ & $(0,0)$ & $0$ & $1$ \\
$(0, 1, 1)$ & $(0,4)$ & $11010, 3734$ & $5$ & $(1, 0, -1)$ & $(0,3)$ & $11602, 1960$ & $4$ \\
$(0, 2, 0)$ & $(0,2)$ & $9946, 1348$ & $2$ & $(1, 1, -1)$ & $(0,1)$ & $5694, 0$ & $2$ \\
\end{longtable}

\endgroup

\begingroup
\tablefont
\setlength{\tabcolsep}{2pt}
\renewcommand{\arraystretch}{1}
\setlength{\LTpre}{2pt}
\setlength{\LTpost}{3pt}
\setlength{\LTcapwidth}{\linewidth}
\tablefont

\begin{longtable}{@{}>{\centering\arraybackslash}p{\dimexpr 0.135000\linewidth-2\tabcolsep\relax}>{\centering\arraybackslash}p{\dimexpr 0.105000\linewidth-2\tabcolsep\relax}>{\centering\arraybackslash}p{\dimexpr 0.190000\linewidth-2\tabcolsep\relax}>{\centering\arraybackslash}p{\dimexpr 0.070000\linewidth-2\tabcolsep\relax}>{\centering\arraybackslash}p{\dimexpr 0.135000\linewidth-2\tabcolsep\relax}>{\centering\arraybackslash}p{\dimexpr 0.105000\linewidth-2\tabcolsep\relax}>{\centering\arraybackslash}p{\dimexpr 0.190000\linewidth-2\tabcolsep\relax}>{\centering\arraybackslash}p{\dimexpr 0.070000\linewidth-2\tabcolsep\relax}@{}}
\caption[Endpoint checks, stage 15]{Endpoint checks, stage 15. Angles (in degrees) $(60/7, 360/7, 120)$; $K=33$; $(k_1,k_2,k_3)=(1, 6, 14)$, $\Lambda=21$; $L=462$, $L(u,v,w,z)=(14719,975,-2653,9345)$; $N\le 621/22$; 21 candidate vertex types. $\boldsymbol\mu=(2, 2, 1)$. The neighbor-number constraint excludes the following 2 vertex types: [(36, 1, 0, 0, 37), (42, 0, 0, 0, 42)]. $\boldsymbol\mu=(2, 2, 1)$. Additional vertex types deleted in this round: 0; elimination terminates.}\label{tab:22}\\[-2pt]
\toprule
$(\varepsilon,t,P)$ & $(s_-,s_+)$ & $\Delta_-,\Delta_+$ & Types & $(\varepsilon,t,P)$ & $(s_-,s_+)$ & $\Delta_-,\Delta_+$ & Types \\
\midrule
\endfirsthead
\caption[]{Continued}\\[-2pt]
\toprule
$(\varepsilon,t,P)$ & $(s_-,s_+)$ & $\Delta_-,\Delta_+$ & Types & $(\varepsilon,t,P)$ & $(s_-,s_+)$ & $\Delta_-,\Delta_+$ & Types \\
\midrule
\endhead
\midrule
\endfoot
\bottomrule
\endlastfoot
$(0, 0, 0)$ & $(2,6)$ & $0, 38724$ & $3$ & $(0, 2, 1)$ & $(1,1)$ & $77028$ & $1$ \\
$(0, 0, 1)$ & $(3,7)$ & $53676, 0$ & $3$ & $(0, 3, 1)$ & $(0,0)$ & $0$ & $1$ \\
$(0, 1, 1)$ & $(0,4)$ & $57078, 40362$ & $5$ & $(1, 0, -1)$ & $(0,3)$ & $96474, 21567$ & $4$ \\
$(0, 2, 0)$ & $(0,2)$ & $97272, 14994$ & $2$ & $(1, 1, -1)$ & $(0,1)$ & $61929, 0$ & $2$ \\
\end{longtable}

\endgroup

\begingroup
\tablefont
\setlength{\tabcolsep}{2pt}
\renewcommand{\arraystretch}{1}
\setlength{\LTpre}{2pt}
\setlength{\LTpost}{3pt}
\setlength{\LTcapwidth}{\linewidth}
\tablefont

\begin{longtable}{@{}>{\centering\arraybackslash}p{\dimexpr 0.135000\linewidth-2\tabcolsep\relax}>{\centering\arraybackslash}p{\dimexpr 0.105000\linewidth-2\tabcolsep\relax}>{\centering\arraybackslash}p{\dimexpr 0.190000\linewidth-2\tabcolsep\relax}>{\centering\arraybackslash}p{\dimexpr 0.070000\linewidth-2\tabcolsep\relax}>{\centering\arraybackslash}p{\dimexpr 0.135000\linewidth-2\tabcolsep\relax}>{\centering\arraybackslash}p{\dimexpr 0.105000\linewidth-2\tabcolsep\relax}>{\centering\arraybackslash}p{\dimexpr 0.190000\linewidth-2\tabcolsep\relax}>{\centering\arraybackslash}p{\dimexpr 0.070000\linewidth-2\tabcolsep\relax}@{}}
\caption[Endpoint checks, stage 16]{Endpoint checks, stage 16. Angles (in degrees) $(60/7, 360/7, 120)$; $K=28$; $(k_1,k_2,k_3)=(1, 6, 14)$, $\Lambda=21$; $L=21$, $L(u,v,w,z)=(569,54,-98,357)$; $N\le 25$; 19 candidate vertex types. $\boldsymbol\mu=(2, 2, 1)$. The neighbor-number constraint excludes the following 4 vertex types: [(28, 0, 1, 0, 29), (30, 2, 0, 0, 32), (36, 1, 0, 0, 37), (42, 0, 0, 0, 42)]. $\boldsymbol\mu=(2, 2, 1)$. Additional vertex types deleted in this round: 0; elimination terminates.}\label{tab:23}\\[-2pt]
\toprule
$(\varepsilon,t,P)$ & $(s_-,s_+)$ & $\Delta_-,\Delta_+$ & Types & $(\varepsilon,t,P)$ & $(s_-,s_+)$ & $\Delta_-,\Delta_+$ & Types \\
\midrule
\endfirsthead
\caption[]{Continued}\\[-2pt]
\toprule
$(\varepsilon,t,P)$ & $(s_-,s_+)$ & $\Delta_-,\Delta_+$ & Types & $(\varepsilon,t,P)$ & $(s_-,s_+)$ & $\Delta_-,\Delta_+$ & Types \\
\midrule
\endhead
\midrule
\endfoot
\bottomrule
\endlastfoot
$(0, 0, 0)$ & $(4,6)$ & $1218, 1218$ & $2$ & $(0, 2, 1)$ & $(1,1)$ & $2688$ & $1$ \\
$(0, 0, 1)$ & $(3,7)$ & $0, 0$ & $3$ & $(0, 3, 1)$ & $(0,0)$ & $0$ & $1$ \\
$(0, 1, 1)$ & $(1,4)$ & $1743, 1428$ & $4$ & $(1, 0, -1)$ & $(0,3)$ & $2352, 777$ & $4$ \\
$(0, 2, 0)$ & $(0,2)$ & $3066, 546$ & $2$ & $(1, 1, -1)$ & $(0,1)$ & $2205, 0$ & $2$ \\
\end{longtable}

\endgroup

\begingroup
\tablefont
\setlength{\tabcolsep}{2pt}
\renewcommand{\arraystretch}{1}
\setlength{\LTpre}{2pt}
\setlength{\LTpost}{3pt}
\setlength{\LTcapwidth}{\linewidth}
\tablefont

\begin{longtable}{@{}>{\centering\arraybackslash}p{\dimexpr 0.135000\linewidth-2\tabcolsep\relax}>{\centering\arraybackslash}p{\dimexpr 0.105000\linewidth-2\tabcolsep\relax}>{\centering\arraybackslash}p{\dimexpr 0.190000\linewidth-2\tabcolsep\relax}>{\centering\arraybackslash}p{\dimexpr 0.070000\linewidth-2\tabcolsep\relax}>{\centering\arraybackslash}p{\dimexpr 0.135000\linewidth-2\tabcolsep\relax}>{\centering\arraybackslash}p{\dimexpr 0.105000\linewidth-2\tabcolsep\relax}>{\centering\arraybackslash}p{\dimexpr 0.190000\linewidth-2\tabcolsep\relax}>{\centering\arraybackslash}p{\dimexpr 0.070000\linewidth-2\tabcolsep\relax}@{}}
\caption[Endpoint checks, stage 17]{Endpoint checks, stage 17. Angles (in degrees) $(60/7, 360/7, 120)$; $K=25$; $(k_1,k_2,k_3)=(1, 6, 14)$, $\Lambda=21$; $L=91$, $L(u,v,w,z)=(2166,263,-357,1344)$; $N\le 296/13$; 17 candidate vertex types. $\boldsymbol\mu=(2, 2, 1)$. The neighbor-number constraint excludes the following 6 vertex types: [(22, 1, 1, 0, 24), (24, 3, 0, 0, 27), (28, 0, 1, 0, 29), (30, 2, 0, 0, 32), (36, 1, 0, 0, 37), (42, 0, 0, 0, 42)]. $\boldsymbol\mu=(2, 2, 1)$. Additional vertex types deleted in this round: 0; elimination terminates.}\label{tab:24}\\[-2pt]
\toprule
$(\varepsilon,t,P)$ & $(s_-,s_+)$ & $\Delta_-,\Delta_+$ & Types & $(\varepsilon,t,P)$ & $(s_-,s_+)$ & $\Delta_-,\Delta_+$ & Types \\
\midrule
\endfirsthead
\caption[]{Continued}\\[-2pt]
\toprule
$(\varepsilon,t,P)$ & $(s_-,s_+)$ & $\Delta_-,\Delta_+$ & Types & $(\varepsilon,t,P)$ & $(s_-,s_+)$ & $\Delta_-,\Delta_+$ & Types \\
\midrule
\endhead
\midrule
\endfoot
\bottomrule
\endlastfoot
$(0, 0, 0)$ & $(4,6)$ & $0, 3654$ & $2$ & $(0, 2, 1)$ & $(1,1)$ & $9212$ & $1$ \\
$(0, 0, 1)$ & $(5,7)$ & $5446, 0$ & $2$ & $(0, 3, 1)$ & $(0,0)$ & $0$ & $1$ \\
$(0, 1, 1)$ & $(2,4)$ & $6776, 4970$ & $3$ & $(1, 0, -1)$ & $(0,3)$ & $4102, 2758$ & $4$ \\
$(0, 2, 0)$ & $(0,2)$ & $9226, 1960$ & $2$ & $(1, 1, -1)$ & $(0,1)$ & $7728, 0$ & $2$ \\
\end{longtable}

\endgroup

\begingroup
\tablefont
\setlength{\tabcolsep}{2pt}
\renewcommand{\arraystretch}{1}
\setlength{\LTpre}{2pt}
\setlength{\LTpost}{3pt}
\setlength{\LTcapwidth}{\linewidth}
\tablefont

\begin{longtable}{@{}>{\centering\arraybackslash}p{\dimexpr 0.135000\linewidth-2\tabcolsep\relax}>{\centering\arraybackslash}p{\dimexpr 0.105000\linewidth-2\tabcolsep\relax}>{\centering\arraybackslash}p{\dimexpr 0.190000\linewidth-2\tabcolsep\relax}>{\centering\arraybackslash}p{\dimexpr 0.070000\linewidth-2\tabcolsep\relax}>{\centering\arraybackslash}p{\dimexpr 0.135000\linewidth-2\tabcolsep\relax}>{\centering\arraybackslash}p{\dimexpr 0.105000\linewidth-2\tabcolsep\relax}>{\centering\arraybackslash}p{\dimexpr 0.190000\linewidth-2\tabcolsep\relax}>{\centering\arraybackslash}p{\dimexpr 0.070000\linewidth-2\tabcolsep\relax}@{}}
\caption[Endpoint checks, stage 18]{Endpoint checks, stage 18. Angles (in degrees) $(60/7, 360/7, 120)$; $K=22$; $(k_1,k_2,k_3)=(1, 6, 14)$, $\Lambda=21$; $L=9$, $L(u,v,w,z)=(173,30,-26,105)$; $N\le 59/3$; 15 candidate vertex types. $\boldsymbol\mu=(2, 2, 1)$. The neighbor-number constraint excludes the following 8 vertex types: [(18, 4, 0, 0, 22), (21, 0, 0, 1, 22), (22, 1, 1, 0, 24), (24, 3, 0, 0, 27), (28, 0, 1, 0, 29), (30, 2, 0, 0, 32), (36, 1, 0, 0, 37), (42, 0, 0, 0, 42)]. $\boldsymbol\mu=(2, 2, 1)$. Additional vertex types deleted in this round: 0; elimination terminates.}\label{tab:25}\\[-2pt]
\toprule
$(\varepsilon,t,P)$ & $(s_-,s_+)$ & $\Delta_-,\Delta_+$ & Types & $(\varepsilon,t,P)$ & $(s_-,s_+)$ & $\Delta_-,\Delta_+$ & Types \\
\midrule
\endfirsthead
\caption[]{Continued}\\[-2pt]
\toprule
$(\varepsilon,t,P)$ & $(s_-,s_+)$ & $\Delta_-,\Delta_+$ & Types & $(\varepsilon,t,P)$ & $(s_-,s_+)$ & $\Delta_-,\Delta_+$ & Types \\
\midrule
\endhead
\midrule
\endfoot
\bottomrule
\endlastfoot
$(0, 0, 0)$ & $(6,6)$ & $138$ & $1$ & $(0, 2, 1)$ & $(1,1)$ & $576$ & $1$ \\
$(0, 0, 1)$ & $(5,7)$ & $36, 0$ & $2$ & $(0, 3, 1)$ & $(0,0)$ & $0$ & $1$ \\
$(0, 1, 1)$ & $(2,4)$ & $0, 324$ & $3$ & $(1, 0, -1)$ & $(1,3)$ & $225, 189$ & $3$ \\
$(0, 2, 0)$ & $(0,2)$ & $354, 138$ & $2$ & $(1, 1, -1)$ & $(0,1)$ & $513, 0$ & $2$ \\
\end{longtable}

\endgroup

\begingroup
\tablefont
\setlength{\tabcolsep}{2pt}
\renewcommand{\arraystretch}{1}
\setlength{\LTpre}{2pt}
\setlength{\LTpost}{3pt}
\setlength{\LTcapwidth}{\linewidth}
\tablefont

\begin{longtable}{@{}>{\centering\arraybackslash}p{\dimexpr 0.135000\linewidth-2\tabcolsep\relax}>{\centering\arraybackslash}p{\dimexpr 0.105000\linewidth-2\tabcolsep\relax}>{\centering\arraybackslash}p{\dimexpr 0.190000\linewidth-2\tabcolsep\relax}>{\centering\arraybackslash}p{\dimexpr 0.070000\linewidth-2\tabcolsep\relax}>{\centering\arraybackslash}p{\dimexpr 0.135000\linewidth-2\tabcolsep\relax}>{\centering\arraybackslash}p{\dimexpr 0.105000\linewidth-2\tabcolsep\relax}>{\centering\arraybackslash}p{\dimexpr 0.190000\linewidth-2\tabcolsep\relax}>{\centering\arraybackslash}p{\dimexpr 0.070000\linewidth-2\tabcolsep\relax}@{}}
\caption[Endpoint checks, stage 19]{Endpoint checks, stage 19. Angles (in degrees) $(10, 50, 120)$; $K=\text{--}$; $(k_1,k_2,k_3)=(1, 5, 12)$, $\Lambda=18$; $L=31$, $L(u,v,w,z)=(1054,-25,-172,609)$; $N\le 857/31$; 23 candidate vertex types.}\label{tab:26}\\[-2pt]
\toprule
$(\varepsilon,t,P)$ & $(s_-,s_+)$ & $\Delta_-,\Delta_+$ & Types & $(\varepsilon,t,P)$ & $(s_-,s_+)$ & $\Delta_-,\Delta_+$ & Types \\
\midrule
\endfirsthead
\caption[]{Continued}\\[-2pt]
\toprule
$(\varepsilon,t,P)$ & $(s_-,s_+)$ & $\Delta_-,\Delta_+$ & Types & $(\varepsilon,t,P)$ & $(s_-,s_+)$ & $\Delta_-,\Delta_+$ & Types \\
\midrule
\endhead
\midrule
\endfoot
\bottomrule
\endlastfoot
$(0, 0, 0)$ & $(0,6)$ & $0, 2454$ & $4$ & $(0, 2, 1)$ & $(1,1)$ & $5138$ & $1$ \\
$(0, 0, 1)$ & $(1,7)$ & $3498, 0$ & $4$ & $(0, 3, 1)$ & $(0,0)$ & $0$ & $1$ \\
$(0, 1, 0)$ & $(1,3)$ & $7460, 4806$ & $2$ & $(1, 0, -1)$ & $(0,3)$ & $8350, 1393$ & $4$ \\
$(0, 1, 1)$ & $(0,4)$ & $7908, 2600$ & $3$ & $(1, 1, -1)$ & $(0,1)$ & $4055, 0$ & $2$ \\
$(0, 2, 0)$ & $(0,2)$ & $7096, 970$ & $2$ &  &  &  &  \\
\end{longtable}

\endgroup

\begingroup
\tablefont
\setlength{\tabcolsep}{2pt}
\renewcommand{\arraystretch}{1}
\setlength{\LTpre}{2pt}
\setlength{\LTpost}{3pt}
\setlength{\LTcapwidth}{\linewidth}
\tablefont

\begin{longtable}{@{}>{\centering\arraybackslash}p{\dimexpr 0.135000\linewidth-2\tabcolsep\relax}>{\centering\arraybackslash}p{\dimexpr 0.105000\linewidth-2\tabcolsep\relax}>{\centering\arraybackslash}p{\dimexpr 0.190000\linewidth-2\tabcolsep\relax}>{\centering\arraybackslash}p{\dimexpr 0.070000\linewidth-2\tabcolsep\relax}>{\centering\arraybackslash}p{\dimexpr 0.135000\linewidth-2\tabcolsep\relax}>{\centering\arraybackslash}p{\dimexpr 0.105000\linewidth-2\tabcolsep\relax}>{\centering\arraybackslash}p{\dimexpr 0.190000\linewidth-2\tabcolsep\relax}>{\centering\arraybackslash}p{\dimexpr 0.070000\linewidth-2\tabcolsep\relax}@{}}
\caption[Endpoint checks, stage 20]{Endpoint checks, stage 20. Angles (in degrees) $(10, 50, 120)$; $K=27$; $(k_1,k_2,k_3)=(1, 5, 12)$, $\Lambda=18$; $L=18$, $L(u,v,w,z)=(437,25,-66,252)$; $N\le 22$; 20 candidate vertex types. $\boldsymbol\mu=(2, 2, 1)$. The neighbor-number constraint excludes the following 3 vertex types: [(26, 2, 0, 0, 28), (31, 1, 0, 0, 32), (36, 0, 0, 0, 36)]. $\boldsymbol\mu=(2, 2, 1)$. Additional vertex types deleted in this round: 0; elimination terminates.}\label{tab:27}\\[-2pt]
\toprule
$(\varepsilon,t,P)$ & $(s_-,s_+)$ & $\Delta_-,\Delta_+$ & Types & $(\varepsilon,t,P)$ & $(s_-,s_+)$ & $\Delta_-,\Delta_+$ & Types \\
\midrule
\endfirsthead
\caption[]{Continued}\\[-2pt]
\toprule
$(\varepsilon,t,P)$ & $(s_-,s_+)$ & $\Delta_-,\Delta_+$ & Types & $(\varepsilon,t,P)$ & $(s_-,s_+)$ & $\Delta_-,\Delta_+$ & Types \\
\midrule
\endhead
\midrule
\endfoot
\bottomrule
\endlastfoot
$(0, 0, 0)$ & $(4,6)$ & $612, 612$ & $2$ & $(0, 2, 1)$ & $(1,1)$ & $1764$ & $1$ \\
$(0, 0, 1)$ & $(3,7)$ & $0, 0$ & $3$ & $(0, 3, 1)$ & $(0,0)$ & $0$ & $1$ \\
$(0, 1, 0)$ & $(1,3)$ & $1080, 1368$ & $2$ & $(1, 0, -1)$ & $(0,3)$ & $1800, 504$ & $4$ \\
$(0, 1, 1)$ & $(0,4)$ & $324, 900$ & $3$ & $(1, 1, -1)$ & $(0,1)$ & $1440, 0$ & $2$ \\
$(0, 2, 0)$ & $(0,2)$ & $2088, 360$ & $2$ &  &  &  &  \\
\end{longtable}

\endgroup

\begingroup
\tablefont
\setlength{\tabcolsep}{2pt}
\renewcommand{\arraystretch}{1}
\setlength{\LTpre}{2pt}
\setlength{\LTpost}{3pt}
\setlength{\LTcapwidth}{\linewidth}
\tablefont

\begin{longtable}{@{}>{\centering\arraybackslash}p{\dimexpr 0.135000\linewidth-2\tabcolsep\relax}>{\centering\arraybackslash}p{\dimexpr 0.105000\linewidth-2\tabcolsep\relax}>{\centering\arraybackslash}p{\dimexpr 0.190000\linewidth-2\tabcolsep\relax}>{\centering\arraybackslash}p{\dimexpr 0.070000\linewidth-2\tabcolsep\relax}>{\centering\arraybackslash}p{\dimexpr 0.135000\linewidth-2\tabcolsep\relax}>{\centering\arraybackslash}p{\dimexpr 0.105000\linewidth-2\tabcolsep\relax}>{\centering\arraybackslash}p{\dimexpr 0.190000\linewidth-2\tabcolsep\relax}>{\centering\arraybackslash}p{\dimexpr 0.070000\linewidth-2\tabcolsep\relax}@{}}
\caption[Endpoint checks, stage 21]{Endpoint checks, stage 21. Angles (in degrees) $(10, 50, 120)$; $K=22$; $(k_1,k_2,k_3)=(1, 5, 12)$, $\Lambda=18$; $L=39$, $L(u,v,w,z)=(859,74,-126,495)$; $N\le 269/13$; 17 candidate vertex types. $\boldsymbol\mu=(2, 2, 1)$. The neighbor-number constraint excludes the following 6 vertex types: [(19, 1, 1, 0, 21), (21, 3, 0, 0, 24), (24, 0, 1, 0, 25), (26, 2, 0, 0, 28), (31, 1, 0, 0, 32), (36, 0, 0, 0, 36)]. $\boldsymbol\mu=(2, 2, 1)$. Additional vertex types deleted in this round: 0; elimination terminates.}\label{tab:28}\\[-2pt]
\toprule
$(\varepsilon,t,P)$ & $(s_-,s_+)$ & $\Delta_-,\Delta_+$ & Types & $(\varepsilon,t,P)$ & $(s_-,s_+)$ & $\Delta_-,\Delta_+$ & Types \\
\midrule
\endfirsthead
\caption[]{Continued}\\[-2pt]
\toprule
$(\varepsilon,t,P)$ & $(s_-,s_+)$ & $\Delta_-,\Delta_+$ & Types & $(\varepsilon,t,P)$ & $(s_-,s_+)$ & $\Delta_-,\Delta_+$ & Types \\
\midrule
\endhead
\midrule
\endfoot
\bottomrule
\endlastfoot
$(0, 0, 0)$ & $(4,6)$ & $0, 918$ & $2$ & $(0, 2, 1)$ & $(1,1)$ & $3210$ & $1$ \\
$(0, 0, 1)$ & $(5,7)$ & $1578, 0$ & $2$ & $(0, 3, 1)$ & $(0,0)$ & $0$ & $1$ \\
$(0, 1, 0)$ & $(3,3)$ & $2250$ & $1$ & $(1, 0, -1)$ & $(0,3)$ & $2370, 939$ & $4$ \\
$(0, 1, 1)$ & $(2,4)$ & $2598, 1644$ & $2$ & $(1, 1, -1)$ & $(0,1)$ & $2661, 0$ & $2$ \\
$(0, 2, 0)$ & $(0,2)$ & $3504, 678$ & $2$ &  &  &  &  \\
\end{longtable}

\endgroup

\begingroup
\tablefont
\setlength{\tabcolsep}{2pt}
\renewcommand{\arraystretch}{1}
\setlength{\LTpre}{2pt}
\setlength{\LTpost}{3pt}
\setlength{\LTcapwidth}{\linewidth}
\tablefont

\begin{longtable}{@{}>{\centering\arraybackslash}p{\dimexpr 0.135000\linewidth-2\tabcolsep\relax}>{\centering\arraybackslash}p{\dimexpr 0.105000\linewidth-2\tabcolsep\relax}>{\centering\arraybackslash}p{\dimexpr 0.190000\linewidth-2\tabcolsep\relax}>{\centering\arraybackslash}p{\dimexpr 0.070000\linewidth-2\tabcolsep\relax}>{\centering\arraybackslash}p{\dimexpr 0.135000\linewidth-2\tabcolsep\relax}>{\centering\arraybackslash}p{\dimexpr 0.105000\linewidth-2\tabcolsep\relax}>{\centering\arraybackslash}p{\dimexpr 0.190000\linewidth-2\tabcolsep\relax}>{\centering\arraybackslash}p{\dimexpr 0.070000\linewidth-2\tabcolsep\relax}@{}}
\caption[Endpoint checks, stage 22]{Endpoint checks, stage 22. Angles (in degrees) $(12, 48, 120)$; $K=\text{--}$; $(k_1,k_2,k_3)=(1, 4, 10)$, $\Lambda=15$; $L=11$, $L(u,v,w,z)=(308,25,-53,192)$; $N\le 280/11$; 24 candidate vertex types.}\label{tab:29}\\[-2pt]
\toprule
$(\varepsilon,t,P)$ & $(s_-,s_+)$ & $\Delta_-,\Delta_+$ & Types & $(\varepsilon,t,P)$ & $(s_-,s_+)$ & $\Delta_-,\Delta_+$ & Types \\
\midrule
\endfirsthead
\caption[]{Continued}\\[-2pt]
\toprule
$(\varepsilon,t,P)$ & $(s_-,s_+)$ & $\Delta_-,\Delta_+$ & Types & $(\varepsilon,t,P)$ & $(s_-,s_+)$ & $\Delta_-,\Delta_+$ & Types \\
\midrule
\endhead
\midrule
\endfoot
\bottomrule
\endlastfoot
$(0, 0, 0)$ & $(0,6)$ & $0, 678$ & $4$ & $(0, 2, 1)$ & $(1,1)$ & $1266$ & $1$ \\
$(0, 0, 1)$ & $(1,7)$ & $800, 290$ & $4$ & $(0, 3, 1)$ & $(0,0)$ & $0$ & $1$ \\
$(0, 1, 1)$ & $(0,5)$ & $1910, 0$ & $6$ & $(1, 0, -1)$ & $(0,3)$ & $1964, 422$ & $4$ \\
$(0, 2, 0)$ & $(0,2)$ & $1654, 296$ & $2$ & $(1, 1, -1)$ & $(0,1)$ & $910, 0$ & $2$ \\
\end{longtable}

\endgroup

\begingroup
\tablefont
\setlength{\tabcolsep}{2pt}
\renewcommand{\arraystretch}{1}
\setlength{\LTpre}{2pt}
\setlength{\LTpost}{3pt}
\setlength{\LTcapwidth}{\linewidth}
\tablefont

\begin{longtable}{@{}>{\centering\arraybackslash}p{\dimexpr 0.135000\linewidth-2\tabcolsep\relax}>{\centering\arraybackslash}p{\dimexpr 0.105000\linewidth-2\tabcolsep\relax}>{\centering\arraybackslash}p{\dimexpr 0.190000\linewidth-2\tabcolsep\relax}>{\centering\arraybackslash}p{\dimexpr 0.070000\linewidth-2\tabcolsep\relax}>{\centering\arraybackslash}p{\dimexpr 0.135000\linewidth-2\tabcolsep\relax}>{\centering\arraybackslash}p{\dimexpr 0.105000\linewidth-2\tabcolsep\relax}>{\centering\arraybackslash}p{\dimexpr 0.190000\linewidth-2\tabcolsep\relax}>{\centering\arraybackslash}p{\dimexpr 0.070000\linewidth-2\tabcolsep\relax}@{}}
\caption[Endpoint checks, stage 23]{Endpoint checks, stage 23. Angles (in degrees) $(12, 48, 120)$; $K=25$; $(k_1,k_2,k_3)=(1, 4, 10)$, $\Lambda=15$; $L=240$, $L(u,v,w,z)=(5702,638,-925,3495)$; $N\le 361/16$; 22 candidate vertex types. $\boldsymbol\mu=(2, 2, 1)$. The neighbor-number constraint excludes the following 2 vertex types: [(26, 1, 0, 0, 27), (30, 0, 0, 0, 30)]. $\boldsymbol\mu=(2, 2, 1)$. Additional vertex types deleted in this round: 0; elimination terminates.}\label{tab:30}\\[-2pt]
\toprule
$(\varepsilon,t,P)$ & $(s_-,s_+)$ & $\Delta_-,\Delta_+$ & Types & $(\varepsilon,t,P)$ & $(s_-,s_+)$ & $\Delta_-,\Delta_+$ & Types \\
\midrule
\endfirsthead
\caption[]{Continued}\\[-2pt]
\toprule
$(\varepsilon,t,P)$ & $(s_-,s_+)$ & $\Delta_-,\Delta_+$ & Types & $(\varepsilon,t,P)$ & $(s_-,s_+)$ & $\Delta_-,\Delta_+$ & Types \\
\midrule
\endhead
\midrule
\endfoot
\bottomrule
\endlastfoot
$(0, 0, 0)$ & $(2,6)$ & $0, 9240$ & $3$ & $(0, 2, 1)$ & $(1,1)$ & $21375$ & $1$ \\
$(0, 0, 1)$ & $(3,7)$ & $12285, 4245$ & $3$ & $(0, 3, 1)$ & $(0,0)$ & $0$ & $1$ \\
$(0, 1, 1)$ & $(0,5)$ & $20850, 0$ & $6$ & $(1, 0, -1)$ & $(0,3)$ & $28275, 7125$ & $4$ \\
$(0, 2, 0)$ & $(0,2)$ & $26370, 5070$ & $2$ & $(1, 1, -1)$ & $(0,1)$ & $15690, 0$ & $2$ \\
\end{longtable}

\endgroup

\begingroup
\tablefont
\setlength{\tabcolsep}{2pt}
\renewcommand{\arraystretch}{1}
\setlength{\LTpre}{2pt}
\setlength{\LTpost}{3pt}
\setlength{\LTcapwidth}{\linewidth}
\tablefont

\begin{longtable}{@{}>{\centering\arraybackslash}p{\dimexpr 0.135000\linewidth-2\tabcolsep\relax}>{\centering\arraybackslash}p{\dimexpr 0.105000\linewidth-2\tabcolsep\relax}>{\centering\arraybackslash}p{\dimexpr 0.190000\linewidth-2\tabcolsep\relax}>{\centering\arraybackslash}p{\dimexpr 0.070000\linewidth-2\tabcolsep\relax}>{\centering\arraybackslash}p{\dimexpr 0.135000\linewidth-2\tabcolsep\relax}>{\centering\arraybackslash}p{\dimexpr 0.105000\linewidth-2\tabcolsep\relax}>{\centering\arraybackslash}p{\dimexpr 0.190000\linewidth-2\tabcolsep\relax}>{\centering\arraybackslash}p{\dimexpr 0.070000\linewidth-2\tabcolsep\relax}@{}}
\caption[Endpoint checks, stage 24]{Endpoint checks, stage 24. Angles (in degrees) $(12, 48, 120)$; $K=22$; $(k_1,k_2,k_3)=(1, 4, 10)$, $\Lambda=15$; $L=45$, $L(u,v,w,z)=(944,131,-145,570)$; $N\le 62/3$; 20 candidate vertex types. $\boldsymbol\mu=(2, 2, 1)$. The neighbor-number constraint excludes the following 4 vertex types: [(20, 0, 1, 0, 21), (22, 2, 0, 0, 24), (26, 1, 0, 0, 27), (30, 0, 0, 0, 30)]. $\boldsymbol\mu=(2, 2, 1)$. Additional vertex types deleted in this round: 0; elimination terminates.}\label{tab:31}\\[-2pt]
\toprule
$(\varepsilon,t,P)$ & $(s_-,s_+)$ & $\Delta_-,\Delta_+$ & Types & $(\varepsilon,t,P)$ & $(s_-,s_+)$ & $\Delta_-,\Delta_+$ & Types \\
\midrule
\endfirsthead
\caption[]{Continued}\\[-2pt]
\toprule
$(\varepsilon,t,P)$ & $(s_-,s_+)$ & $\Delta_-,\Delta_+$ & Types & $(\varepsilon,t,P)$ & $(s_-,s_+)$ & $\Delta_-,\Delta_+$ & Types \\
\midrule
\endhead
\midrule
\endfoot
\bottomrule
\endlastfoot
$(0, 0, 0)$ & $(4,6)$ & $780, 1050$ & $2$ & $(0, 2, 1)$ & $(1,1)$ & $3240$ & $1$ \\
$(0, 0, 1)$ & $(3,7)$ & $0, 540$ & $3$ & $(0, 3, 1)$ & $(0,0)$ & $0$ & $1$ \\
$(0, 1, 1)$ & $(1,5)$ & $2700, 0$ & $5$ & $(1, 0, -1)$ & $(0,3)$ & $3510, 1080$ & $4$ \\
$(0, 2, 0)$ & $(0,2)$ & $3750, 780$ & $2$ & $(1, 1, -1)$ & $(0,1)$ & $2430, 0$ & $2$ \\
\end{longtable}

\endgroup

\begingroup
\tablefont
\setlength{\tabcolsep}{2pt}
\renewcommand{\arraystretch}{1}
\setlength{\LTpre}{2pt}
\setlength{\LTpost}{3pt}
\setlength{\LTcapwidth}{\linewidth}
\tablefont

\begin{longtable}{@{}>{\centering\arraybackslash}p{\dimexpr 0.135000\linewidth-2\tabcolsep\relax}>{\centering\arraybackslash}p{\dimexpr 0.105000\linewidth-2\tabcolsep\relax}>{\centering\arraybackslash}p{\dimexpr 0.190000\linewidth-2\tabcolsep\relax}>{\centering\arraybackslash}p{\dimexpr 0.070000\linewidth-2\tabcolsep\relax}>{\centering\arraybackslash}p{\dimexpr 0.135000\linewidth-2\tabcolsep\relax}>{\centering\arraybackslash}p{\dimexpr 0.105000\linewidth-2\tabcolsep\relax}>{\centering\arraybackslash}p{\dimexpr 0.190000\linewidth-2\tabcolsep\relax}>{\centering\arraybackslash}p{\dimexpr 0.070000\linewidth-2\tabcolsep\relax}@{}}
\caption[Endpoint checks, stage 25]{Endpoint checks, stage 25. Angles (in degrees) $(120/7, 300/7, 120)$; $K=\text{--}$; $(k_1,k_2,k_3)=(2, 5, 14)$, $\Lambda=21$; $L=9$, $L(u,v,w,z)=(166,37,-26,105)$; $N\le 59/3$; 14 candidate vertex types.}\label{tab:32}\\[-2pt]
\toprule
$(\varepsilon,t,P)$ & $(s_-,s_+)$ & $\Delta_-,\Delta_+$ & Types & $(\varepsilon,t,P)$ & $(s_-,s_+)$ & $\Delta_-,\Delta_+$ & Types \\
\midrule
\endfirsthead
\caption[]{Continued}\\[-2pt]
\toprule
$(\varepsilon,t,P)$ & $(s_-,s_+)$ & $\Delta_-,\Delta_+$ & Types & $(\varepsilon,t,P)$ & $(s_-,s_+)$ & $\Delta_-,\Delta_+$ & Types \\
\midrule
\endhead
\midrule
\endfoot
\bottomrule
\endlastfoot
$(0, 0, 0)$ & $(2,6)$ & $138, 138$ & $2$ & $(0, 2, 1)$ & $(0,0)$ & $648$ & $1$ \\
$(0, 0, 1)$ & $(0,8)$ & $0, 0$ & $3$ & $(0, 3, 1)$ & $(0,0)$ & $0$ & $1$ \\
$(0, 1, 1)$ & $(0,4)$ & $648, 324$ & $3$ & $(1, 0, -1)$ & $(1,3)$ & $540, 189$ & $2$ \\
$(0, 2, 0)$ & $(2,2)$ & $138$ & $1$ & $(1, 1, -1)$ & $(1,1)$ & $0$ & $1$ \\
\end{longtable}%

\endgroup

\phantomsection\label{app:special-slacks}

Table~\ref{tab:39} lists all candidate vertex types and inequality slacks for the special right triangle.
Using the vertex types $\tau_B,\tau_C,\tau_R$ from Proposition~\ref{prop:18}, define
\[
\begin{aligned}
\Delta_0&=90r+14s+10t-5\varepsilon+5\cdot\mathbf1_{\tau_B}-5q(q-2),\\
\delta&=18r+3s+2t-q(q-2),\\
\Delta_1&=8\delta-9s-29t-29\varepsilon.
\end{aligned}
\]
The deficits of vertex types $(20,0,0;0)$, $\tau_B$, $\tau_C$ and $\tau_R$ are counted separately in the proof. The $\Delta_1\ (\delta)$ columns of Table~\ref{tab:39} give $\delta$ for these four types and the slack $\Delta_1$ of inequality~\eqref{eq:deficit-estimate} for all others; the $\Delta_0$ columns give the first-round inequality slacks for every candidate vertex type.

\begingroup
\tablefont
\setlength{\tabcolsep}{2pt}
\renewcommand{\arraystretch}{1}
\setlength{\LTpre}{2pt}
\setlength{\LTpost}{3pt}
\setlength{\LTcapwidth}{\linewidth}
\tablefont

\begin{longtable}{@{}>{\centering\arraybackslash}p{\dimexpr 0.210000\linewidth-2\tabcolsep\relax}>{\centering\arraybackslash}p{\dimexpr 0.060000\linewidth-2\tabcolsep\relax}>{\centering\arraybackslash}p{\dimexpr 0.100000\linewidth-2\tabcolsep\relax}>{\centering\arraybackslash}p{\dimexpr 0.130000\linewidth-2\tabcolsep\relax}>{\centering\arraybackslash}p{\dimexpr 0.210000\linewidth-2\tabcolsep\relax}>{\centering\arraybackslash}p{\dimexpr 0.060000\linewidth-2\tabcolsep\relax}>{\centering\arraybackslash}p{\dimexpr 0.100000\linewidth-2\tabcolsep\relax}>{\centering\arraybackslash}p{\dimexpr 0.130000\linewidth-2\tabcolsep\relax}@{}}
\caption{Complete integer slacks for the exceptional right triangle}\label{tab:39}\\[-2pt]
\toprule
Vertex type & $q$ & $\Delta_0$ & $\Delta_1\ (\delta)$ & Vertex type & $q$ & $\Delta_0$ & $\Delta_1\ (\delta)$ \\
\midrule
\endfirsthead
\caption[]{Continued}\\[-2pt]
\toprule
Vertex type & $q$ & $\Delta_0$ & $\Delta_1\ (\delta)$ & Vertex type & $q$ & $\Delta_0$ & $\Delta_1\ (\delta)$ \\
\midrule
\endhead
\midrule
\endfoot
\bottomrule
\endlastfoot
$(20,0,0;0)$ & $20$ & $0$ & $\delta=0$ & $(6,1,2;0)$ & $9$ & $259$ & $349$ \\
$(16,1,0;0)$ & $17$ & $179$ & $279$ & $(2,2,2;0)$ & $6$ & $108$ & $100$ \\
$(12,2,0;0)$ & $14$ & $268$ & $414$ & $(5,0,3;0)$ & $8$ & $240$ & $297$ \\
$(8,3,0;0)$ & $11$ & $267$ & $405$ & $(1,1,3;0)$ & $5$ & $59$ & $0$ \\
$(4,4,0;0)$ & $8$ & $176$ & $252$ & $(0,0,4;0)$ & $4$ & $0$ & $\delta=0$ \\
$(0,5,0;0)$ & $5$ & $0$ & $\delta=0$ & $(10,0,0;1)$ & $11$ & $400$ & $619$ \\
$(15,0,1;0)$ & $16$ & $240$ & $355$ & $(6,1,0;1)$ & $8$ & $309$ & $466$ \\
$(11,1,1;0)$ & $13$ & $299$ & $442$ & $(2,2,0;1)$ & $5$ & $128$ & $169$ \\
$(7,2,1;0)$ & $10$ & $268$ & $385$ & $(5,0,1;1)$ & $7$ & $280$ & $398$ \\
$(3,3,1;0)$ & $7$ & $147$ & $184$ & $(1,1,1;1)$ & $4$ & $69$ & $53$ \\
$(10,0,2;0)$ & $12$ & $320$ & $454$ & $(0,0,2;1)$ & $3$ & $0$ & $\delta=1$ \\
\end{longtable}

\endgroup

\phantomsection\label{app:rigid-cert}

For the ordinary weighted inequalities in the equality classification of Section~\ref{sec:rigidity}, define the slack by
\[
\Delta=ur+vs+wt+zP-q(q-2).
\]
Vertex types are written as $(r,s,t;\varepsilon)$, with $q=r+s+t+\varepsilon$ and $P$ as defined in the main text. In the ordinary Cases C1--C22, the four coefficients are taken from the corresponding row of Table~\ref{article-rigid-table-1}; in the special Cases C23 and C24, the slacks are defined by the formulas in the corresponding cases of Table~\ref{app:rigid-slacks-combined}.
The initial list is determined by the vertex angle equation and the bounds on $r,s$ in the main text. At each round, recompute $\mu_i$ on the current list and simultaneously delete the vertex types violating the constraint of $21$ neighbors on each tile, continuing until the list stabilizes. Table~\ref{app:rigid-slacks-combined} gives, for each case, the minima and deletions at each round, followed by the slacks of the final candidate vertex types. These candidate vertex types together with all deleted vertex types form the initial list. 
The table lists $378$ slacks for the ordinary weighted inequalities and $44$ for the special inequalities, all nonnegative. 

\begingroup
\tablefont
\setlength{\tabcolsep}{2pt}
\renewcommand{\arraystretch}{1}
\setlength{\LTpre}{2pt}
\setlength{\LTpost}{3pt}
\setlength{\LTcapwidth}{\linewidth}
\tablefont

\renewcommand{\arraystretch}{0.95}


\endgroup
\end{document}